\documentclass[a4paper,10pt]{article}
\usepackage{mathtext}
\usepackage[T1,T2A]{fontenc}
\usepackage[cp1251]{inputenc}
\usepackage[english]{babel}
\usepackage{amsmath}
\usepackage{amsfonts}
\usepackage{amssymb}
\usepackage{mathrsfs}
\usepackage{amsthm}
\usepackage{titling}
\usepackage{multicol}
\usepackage{graphicx}
\usepackage{wrapfig}
\usepackage{appendix}

\usepackage{euscript}

\DeclareMathOperator{\interior}{int}

\DeclareMathOperator{\dist}{dist}

\newcommand{\df}{\buildrel\mathrm{def}\over=}
\newcommand{\Bell}{\boldsymbol{B}}
\newcommand{\Bellb}{\boldsymbol{B}^{\mathrm{b}}}
\newcommand{\Class}{\boldsymbol{A}}
\newcommand{\Classd}{\boldsymbol{A}^{\mathrm{dyad}}}
\newcommand{\Classb}{\boldsymbol{A}^{\mathrm{b}}}
\newcommand{\BMO}{\mathrm{BMO}}
\newcommand{\eps}{\varepsilon}
\DeclareMathOperator{\cl}{cl}
\DeclareMathOperator{\pr}{\mathrm{pr}}
\DeclareMathOperator{\conv}{conv}
\DeclareMathOperator{\E}{\mathbb{E}}
\newcommand{\av}[2]{\langle {#1}\rangle_{{}_{#2}}}

\renewcommand{\leq}{\leqslant}
\renewcommand{\geq}{\geqslant}
\renewcommand{\emptyset}{\varnothing}
\newcommand{\FixedBoundary}{\partial_{\mathrm{fixed}}}
\newcommand{\FreeBoundary}{\partial_{\mathrm{free}}}

\newcommand{\ul}{u_{\mathrm{L}}}
\newcommand{\ur}{u_{\mathrm{R}}}

\newcommand{\BG}{\mathfrak{B}}
\newcommand{\BM}{\EuScript{B}}
\newcommand{\BMb}{\EuScript{B}^{\mathrm{b}}}

\newcommand{\M}{\EuScript{M}}
\newcommand{\Mb}{\EuScript{M}^{\mathrm{b}}}

\newtheorem*{MainT}{Theorem}
\newtheorem{Le}{Lemma}[section]
\newtheorem{Def}[Le]{Definition}
\newtheorem{St}[Le]{Proposition}
\newtheorem{Th}[Le]{Theorem}
\newtheorem{Cor}[Le]{Corollary}
\newtheorem{Rem}[Le]{Remark}
\newtheorem{Fact}[Le]{Fact}
\newtheorem{Conj}[Le]{Conjecture}

\numberwithin{equation}{section}

\begin{document}
\author{D.~M.~Stolyarov\thanks{The research is supported by the grant of the Russian Science Foundation (project 14-21-00035).} \and P.~B.~Zatitskiy\thanksmark{1}}
\title{Theory of locally concave functions and its applications to sharp estimates of integral functionals}

\maketitle
\begin{abstract}
We prove a duality theorem the computation of certain Bellman functions is usually based on. As a byproduct, we obtain sharp results about the norms of monotonic rearrangements. The main novelty of our approach is a special class of martingales and an extremal problem on this class, which is dual to the minimization problem for locally concave functions.  
\end{abstract}

\tableofcontents

\section{Introduction}\label{s1}

\subsection{Setting and main ideas}\label{s11}
Some results of the present paper have been announced in~\cite{IOSVZ2}. However, the setting we use here is a bit different, so we give it in full detial.

\begin{figure}[h!]
\includegraphics[height=9.5cm]{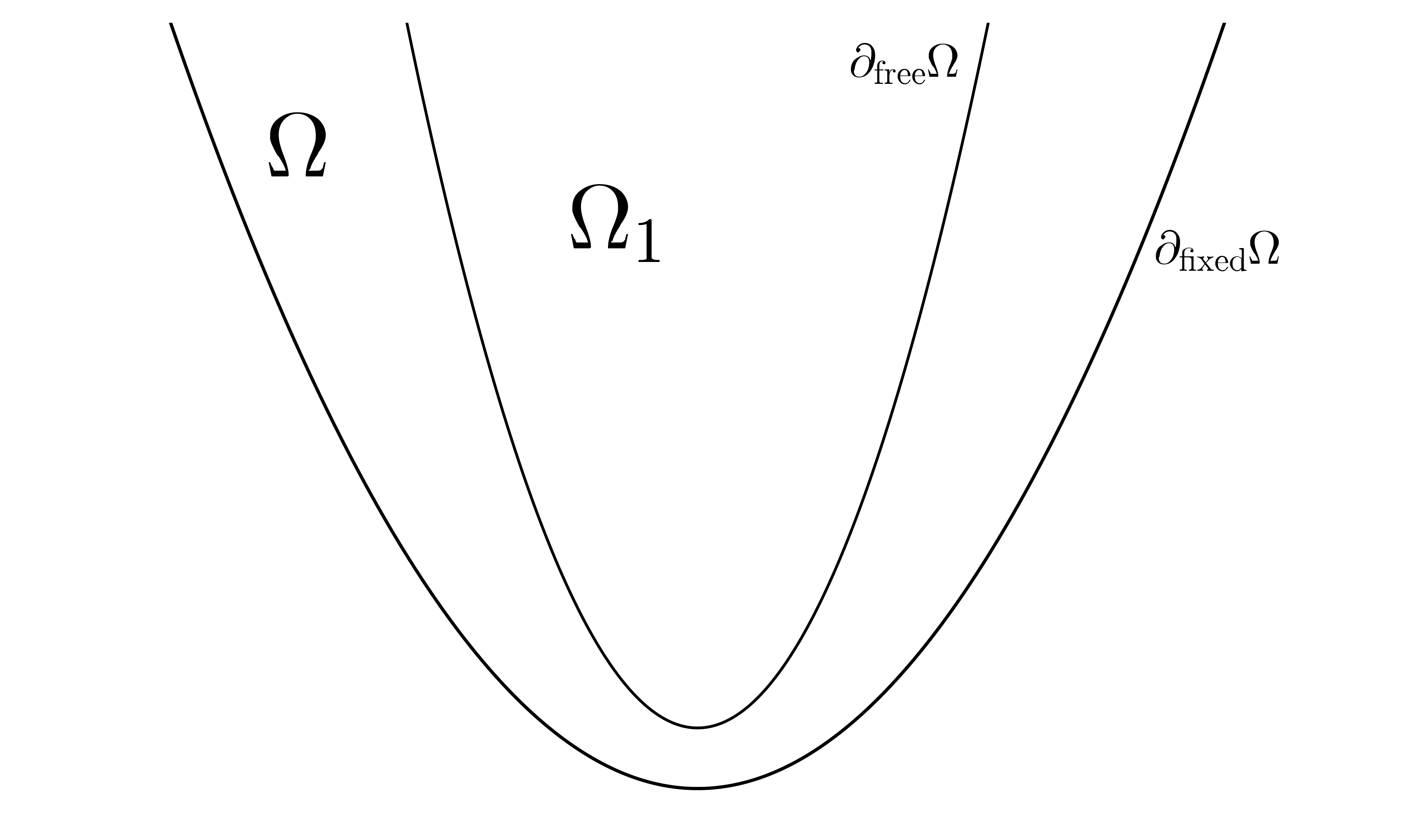}
\caption{Set~$\Omega$ and its boundaries.}
\label{fig:strip}
\end{figure}

Let~$\Omega_0$ be a non-empty open convex subset of~$\mathbb{R}^2$ that does not contain lines. Let~$\Omega_1$ be another open convex subset of~$\mathbb{R}^2$ such that~$\cl \Omega_1 \subset \Omega_0$; the symbol~$\cl$ denotes the closure. We define the domain~$\Omega$ as~$\cl (\Omega_0 \setminus  \Omega_1)$ (the word ``domain'' comes from ``domain of a function''), see Figure~\ref{fig:strip} for visualization, and the class of summable~$\mathbb{R}^2$-valued functions on an interval~$I \subset \mathbb{R}$:
\begin{equation}\label{AnalyticClass}
\Class_{\Omega} = \big\{\varphi \in L^1(I,\mathbb{R}^2)\,\big|\,\, \varphi(I) \subset \partial \Omega_0, \,\,\forall J \hbox{---subinterval of}\,\,I\quad \av{\varphi}{J} \notin \Omega_1\big\}.
\end{equation}
Here~$\av{\varphi}{J} = \frac{1}{|J|}\int_J \varphi(s)\,ds$ is the average of~$\varphi$ over~$J$.
In Subsection~\ref{s13}, we show how the~$\BMO$, the Muckenhoupt classes, and the Gehring classes can be represented in the form~\eqref{AnalyticClass}. Let~$f$ be a bounded from below Borel measurable locally bounded function on~$\partial \Omega_0$\footnote{The function~$f$ is always assumed to be Borel measurable and locally bounded. Sometimes we assume that it is bounded from below, sometimes not. The second case is more interesting from the theoretical point of view. For the first reading, we recommend to assume everywhere that~$f$ is bounded from below.}. We are interested in sharp bounds for the expressions of the form
\begin{equation*}
\av{f(\varphi)}{I}, \quad \varphi \in \Class_{\Omega}.
\end{equation*}
Again, in Subsection~\ref{s13}, we explain how the John--Nirenberg inequality or other inequalities of harmonic analysis can be rewritten as estimations of such an expression. The said estimates are delivered by the corresponding Bellman function,
\begin{equation}\label{BellmanFunction}
\Bell_{\Omega,f}(x) = \sup\big\{\av{f(\varphi)}{I}\,\big|\,\, \av{\varphi}{I} = x,\,\,\varphi \in \Class_{\Omega}\big\}, \quad x \in \Omega.
\end{equation}
The aim of this paper is to prove that this function enjoys good analytic properties. 
\begin{Def}\label{LocalConcavity}
Let~$\omega$ be a subset of~$\mathbb{R}^d$. We call a function~$G\colon w \to \mathbb{R}\cup\{+\infty\}$ locally concave on~$\omega$ provided for every segment~$\ell \subset \omega$ the restriction~$G\big|_{\ell}$ is concave.
\end{Def}
Define the class of functions on~$\Omega$:
\begin{equation}\label{Lambdaclass}
\Lambda_{\Omega,f} = \Big\{G\colon \Omega\to \mathbb{R}\cup\{+\infty\}
\,\Big|\,\, 
G\hbox{ is locally concave on~$\Omega$,} \quad \forall x \in \partial \Omega_0 \quad G(x) \geq f(x)\Big\}.
\end{equation}
The function~$\BG_{\Omega,f}$ is given as follows,
\begin{equation}\label{MinimalLocallyConcave}
\BG_{\Omega,f}(x) = \inf_{G \in \Lambda_{\Omega,f}} G(x), \quad x \in \Omega.
\end{equation}
The main theorem says that~$\Bell = \BG$ (here and in what follows we omit indices of these functions if this does not cause ambiguity). We impose several technical conditions on~$\Omega$ (we will clarify the meaning of the third condition in Subsection~\ref{s31}).
\begin{align}
\label{FirstCondition} &1.\hbox{ Neither~$\partial \Omega_0$ nor~$\partial \Omega_1$ contains linear segments.}\\
\label{SecondCondition} &2.\hbox{ The boundary of~$\Omega_1$ is~$C^2$-smooth.}\\
\label{ThirdCondition} &3.\hbox{ The domain~$\Omega_0$ is unbounded and every ray inside it can be translated to belong to~$\Omega_1$ entirely.}
\end{align}
\begin{MainT}\label{MT}
If~$\Omega$ satisfies technical assumptions \textup{\eqref{FirstCondition}, \eqref{SecondCondition}, \eqref{ThirdCondition},} then~$\Bell = \BG$.
\end{MainT}

The motivation for the problem is given in Subsections~\ref{s12},~\ref{s13}, now we sketch the idea of the proof and provide the structure of the paper. 

The main idea is to introduce the third function that is based on some optimization. Namely, we consider all~$\mathbb{R}^2$-valued martingales that walk inside~$\Omega$ and end their way on~$\partial \Omega_0$. Then, for any point~$x$ in~$\Omega$ we maximize the value~$\E f(M_{\infty})$ over all the martingales starting from~$x$, and denote this supremum by~$\BM(x)$. It is not difficult to prove that the achieved function coincides with~$\BG$, moreover this relation holds true for a general ``good enough'' domain in any dimension\footnote{The authors have a strong belief that this duality may be useful outside the Bellman function theory.}. This construction is described in Section~\ref{s2}.

To prove the main theorem, we prove two inequalities,~$\BG \leq \Bell$ and~$\Bell \leq \BG$. The first one uses ideas of Section~\ref{s2}. Namely, every martingale in~$\Omega$ gives rise to a function belonging to the class~$\Class_{\Omega}$ with the same distribution. Therefore,~$\BM \leq \Bell$, and by results of~Subsection~\ref{s2}, the first inequality follows. To prove the second one, we establish a reverse embedding. It turns out that each function belonging to~$\Class_{\Omega}$ gives rise to a martingale that lives in an extension of~$\Omega$ (i.e. a similar-built domain with strictly smaller~$\Omega_1$), moreover, the difference between this extension and~$\Omega$ can be arbitrary small. Therefore, the inequality~$\Bell \leq \BG$ is almost proved, we have~$\Bell_{\Omega,f} \leq \BG_{\tilde{\Omega},f}$ for any extension~$\tilde{\Omega}$. All this material constitutes Section~\ref{s3}.

To finish the proof of the main theorem, we establish that~$\inf_{\tilde{\Omega}\supset \Omega}\BG_{\tilde{\Omega}} = \BG_{\Omega}$ (the infimum is taken over all the extensions of~$\Omega$). This is done in Section~\ref{s4} for the case where~$\partial \Omega_0$ is~$C^2$-smooth and~$f$ is continuous (Theorem~\ref{SmoothMainTheorem}). 

Section~\ref{s5} treats the case of non-smooth boundary~$\partial \Omega_0$ and non-smooth~$f$. More or less, the result is derived in a classical way: we do some sort of smoothing, apply the already known results for smoothed functions, and then pass to the limit. However, the non-linearity of the problem makes the smoothing non-standard, some geometric tricks are used here. Corollary~\ref{DiscontinuousFunctionNonSmoothBoundary} finishes the proof of the main theorem in full generality.

In Section~\ref{s6}, we give some information for the case where~$f$ is not bounded from below. Now the problem should be re-stated, because \emph{a priori} the value~$\av{f(\varphi)}{I}$ is not well defined (the function~$f(\varphi)$ might be not integrable). The fact that the integral of~$f(\varphi)$ is well defined for all~$\varphi$ belonging to the class~$\Class_{\Omega}$ is equivalent to the finiteness of~$\BG_{\Omega,f_+}$ ($f_+$ stands for the positive part of~$f$). However, the condition that~$\BG_{\Omega,f}$ is finite is not sufficient here. We also give a sufficient summability condition for~$\BG_{\Omega,f}$ to be finite in terms of some maximal function of~$f$  (this condition can be easily verified with the function~$f$ at hand).  

In Section~\ref{s7}, we formulate several conjectures.

There are also three Appendices that collect various supplementary material.

\paragraph{Acknowledgments.} We are grateful to our colleagues Paata Ivanisvili, Alexander Logunov, and Nikolay Osipov for their criticism, and Leonid Slavin for helpful exposition advice. We also thank Sergey~Vladi\-mirovich~Ivanov who suggested the idea of using the projective transform in this context (see Subsection~\ref{s13} and Appendix~\ref{SPTT}).

The second half of this text was written while the first author was visiting Hausdorff Institute for Mathematics, University of Bonn. He thanks HIM for hospitality.

We are grateful to our teacher Vasily~Ivanovich~Vasyunin for his support, advice, and attention to our work.  

\subsection{Historical remarks}\label{s12}
Application of optimization ideas to analytic problems has a long history. In 1984, Burkholder in his seminal paper~\cite{Burkholder} provided sharp estimates for the norm of a martingale transform in~$L^p$. His method was based on a certain extremal problem of finding a minimal diagonally concave function on a special domain in~$\mathbb{R}^3$. Afterwards, there were many papers where similar technique was used for proving different sharp inequalities for martingales, see the book~\cite{Osekowski}, references therein, and the survey paper~\cite{Osekowski3}.

In the mid ninetieth, Nazarov, Treil, and Volberg introduced optimization principles to harmonic analysis. See~\cite{NTV} for the history of the development and also~\cite{NT} as the historically first exposition. The strength of their method was in building \emph{supersolutions},  which still provided good estimates, i.e. finding not exact Bellman functions. Since then, the method has become a standard tool in analysis, for example, see the lecture notes~\cite{Volberg1,Volberg} or the survey~\cite{NTV}.

However, we mention a much earlier work of Hanner~\cite{Hanner} that presents Beurling's proof of the so-called Hanner's inequalities dating back to 1945. In fact, this method is very Bellman-style: one guesses a special function that proves the desired inequality for him. See~\cite{ISZ} for the explanations where is the Bellman function hidden there.

Around the year of 2002,~Slavin  and  Vasyunin independently found the sharp constants in the John--Nirenberg inequality, see~\cite{Slavin,Vasyunin}, and finally~\cite{SV}. Seemingly, this was the first exact Bellman function for a purely  harmonic-analytic problem (we also mention the paper~\cite{Melas} in this context, which appeared a bit later, but works with a different problem of estimating the~$L^p$ norm of maximal operators). The paper~\cite{Vasyunin2} was the first to appear. Following this route, many authors have managed to prove various inequalities, see the papers~\cite{BR,DW,IOSVZ,ISVZ,LSSVZ,Slavin,Osekowski2,Reznikov,SSV,SV,SV2,Vasyunin,Vasyunin2,Vasyunin3,Vasyunin4,Vasyunin5,VV,VV1}. Though the method turned out to be very useful, there was absolutely no theory that allowed to calculate the Bellman functions ``mechanically''. The first steps towards building such a theory were done in~\cite{SV2}. The problems concerning the~$\BMO$-space, seemingly, form the largest group of the solved problems. In~\cite{IOSVZ}, most of them were unified into a single theory (see also the short report~\cite{IOSVZShortReport}). However, for the treatment in the full generality, see the forthcoming paper~\cite{ISVZ} (still on the space~$\BMO$!). 

All the mentioned papers (including the latest papers~\cite{IOSVZ,ISVZ}), in some sense, employ a miracle. The main theorem of the present paper was usually proved for a specific case in such a fashion: it was assumed to be valid, from these assumptions, the author guessed the ``formula'' for the Bellman function, after that he guessed the optimizers, and finally, he verified the concavity, thus proving that his guess for the Bellman function was right (and also proving the main theorem for his particular case). However, each time this took lots of pages of calculations mixed up with magic guesses. So, the main theorem itself was a miracle. Our aim is to find the reasons for it. In some sense, our explanations show that there are no ``harmonic analytic'' Bellman functions, they are hidden Bellman functions for optimization of stochastic processes. On the formal level, we believe that our studies may give additional information about the links between the Burkholder method (as presented in the book~\cite{Osekowski}) and the Bellman function method.

As a byproduct, we obtain sharp inequalities for non-increasing rearrangements of functions in the class considered (in particular, in the~$\BMO$, the Muckenhoupt classes, the Gehring classes), namely, we prove that the non-increasing rearrangement lies in the same class as the function does (so it does not increase the~$\BMO$ norm, or the Muckenhoupt constant, or the Gehring constant). These results are presented in Subsection~\ref{s34}. The first sharp inequality of such kind (that the non-increasing rearrangement does not increase the~$\BMO$-norm of functions) was proved by Klemes,~\cite{Klemes}, and then generalized in~\cite{Korenovskii}. See also~\cite{KS} for a survey on monotonic rearrangements of functions in~$\BMO$. The case of the class~$A_1$ was considered in~\cite{BSW} (and reproved by several authors later). For the classes~$A_p$ and the Gehring (reverse-H\"older) classes, a similar statement was proved in~\cite{Korenovskii2}. Our approach is different from those used in the papers cited above. More or less, results of this type are direct consequences of the constructions lying behind the main theorem.  

\subsection{Particular cases}\label{s13}
\paragraph{The~$\BMO$ space.} We consider the~$\BMO$ space with the quadratic seminorm. Let~$\eps$ be a positive number. Let~$\Omega_0 = \{x\in\mathbb{R}^2\mid x_1^2 < x_2\}$, let~$\Omega_{1} = \{x\in\mathbb{R}^2\mid  x_1^2 + \eps^2 < x_2\}$. The function~$\varphi = (\varphi_1,\varphi_2)\colon I \to \partial\Omega_0$ belongs to the class~$\Class_{\Omega}$ if and only if~$\varphi_1$ (its first coordinate) belongs to~$\BMO_{\eps}$ (the ball of the~$\BMO$ space of radius~$\eps$). Indeed, for any~$t \in I$ we have~$\varphi_2(t) = \varphi_1^2(t)$, therefore, the condition~$\av{\varphi}{J}\notin\Omega_1$ can be rewritten as
\begin{equation*}
\av{\varphi_1^2}{J} \leq \av{\varphi_1}{J}^2 + \eps^2,
\end{equation*}
which is the same as
\begin{equation}\label{AlmostBMO}
\av{\big(\varphi_1 - \av{\varphi_1}{J}\big)^2}{J} \leq \eps^2.
\end{equation}
Now we see that the class~$\Class_{\Omega}$ corresponds to~$\BMO_{\eps}$. The Bellman function~\eqref{BellmanFunction} estimates the functional~$\av{\tilde{f}(\varphi_1)}{I}$, where~$f(\varphi) = \tilde{f}(\varphi_1)$. We address the reader to the paper~\cite{IOSVZ}, where it was explained how do the sharp estimates of these functionals lead to various forms of the John--Nirenberg inequality and the equivalence of the~$L^p$-type norms defining the~$\BMO$. This case is the subject of study for papers~\cite{IOSVZ, ISVZ, LSSVZ, Osekowski2, SV, SV2, Vasyunin4, Vasyunin5}.

\paragraph{Classes~$A_{p_1,p_2}$.} Let~$p_1$ and~$p_2$ be real numbers and let~$Q \geq 1$. Define the domain~$\Omega$ by the formula
\begin{equation*}
\Omega_0 = \{x \in\mathbb{R}^2\mid  x_1,x_2 > 0,\,\, x_2^{\frac{1}{p_2}} < x_1^{\frac{1}{p_1}}\} \quad \hbox{and}\quad \Omega_1 = \{x \in\mathbb{R}^2\mid  x_1,x_2 > 0,\,\, Q x_2^{\frac{1}{p_2}} < x_1^{\frac{1}{p_1}}\}. 
\end{equation*}
We warn the reader that these domains may not satisfy our conditions (if~$p_1p_2 > 0$), however, we will show how to deal with this problem in Appendix~\ref{SPTT} (the idea is to make a projective transform). If the function~$\varphi$ belongs to the class~$\Class_{\Omega}$, then its first coordinate,~$\varphi_1$, belongs to the so-called~$A_{p_1,p_2}$ class. The ``norm'' in this class is defined as
\begin{equation}\label{Apnorm}
[\psi]_{A_{p_1,p_2}} = \sup\limits_{J \subset I} \,\,\,\av{\psi^{p_1}}{J}^{\frac{1}{p_1}}\av{\psi^{p_2}}{J}^{-\frac{1}{p_2}},
\end{equation}
where the supremum is taken over all the subintervals of~$I$.
These classes were introduced in~\cite{Vasyunin3}. If~$p \in (1,\infty)$, then~$A_{1,-\frac{1}{p-1}} = A_p$, where~$A_p$ stands for the classical Muckenhoupt class. The limiting cases~$A_1$ and~$A_{\infty}$ also fit this definition (with Hruschev's ``norm'' on~$A_{\infty}$). When~$p_2 = 1$ and~$p_1 > 1$, the class~$A_{p_1,p_2}$ coincides with the so-called Gehring class (see~\cite{Korenovskii2} or~\cite{KS}). One can see that the functions in the Gehring class are exactly those that satisfy the reverse H\"older inequality. Sometimes, the Gehring class is called the reverse-H\"older class. Estimates of integral functionals as provided by the Bellman function~\eqref{BellmanFunction} lead to various sharp forms of the reverse H\"older inequality, see~\cite{Vasyunin3}. These cases were treated in the papers~\cite{BR, DW, Reznikov, Vasyunin2, Vasyunin3}.

\paragraph{Reverse Jensen classes.} These classes were introduced in~\cite{Korenovskii2}. Let~$\Phi\colon \mathbb{R}_+ \to \mathbb{R}_+$ be a convex function. Let~$Q > 1$. Consider the class of functions~$\psi\colon I \to \mathbb{R}_+$ such that
\begin{equation*}
\forall J \subset I \quad \av{\Phi(\psi)}{J} \leq Q \Phi(\av{\psi}{J}).
\end{equation*}
Surely, both Muckenhoupt classes and Gehring classes can be described as certain Reverse Jensen classes. The corresponding domain is~$\{x \in \mathbb{R}^2 \mid x_1,x_2 > 0,\,\, \Phi(x_1) \leq x_2 \leq Q \Phi(x_1)\}$. The case~$\Phi(t) = e^t$ turned out to be useful in the study of the John--Nirenberg inequality for the case of the~$\BMO$~$p$-norm, see~\cite{Slavin3}. Unless~$\lim_{x\to 0+}\Phi(x) =\infty$, the domain does not satisfy the conditions required. We note that these classes have not been studied with the Bellman function method, however, it was proved in~\cite{Korenovskii2} that the monotonic rearrangement does not drop the function out of such a class. In the case where~$\lim_{x\to 0+}\Phi(x) = \lim_{x\to 0+}\Phi'(x) = 0$ and~$\lim_{x \to \infty}\frac{\Phi(x)}{x} = \infty$, the domain satisfies the conditions in Appendix~\ref{SPTT}, thus, the main theorem and the statement about monotonic rearrangement are valid for it. To get the assertion on monotonic rearrangements in the full generality, one may approximate the domain in question by the domains of the type described in Subsection~\ref{s1}.

\section{Martingales on domains}\label{s2}
\subsection{Properties of locally concave functions}\label{s21}
\begin{figure}[h!]
\includegraphics[height=9.5cm]{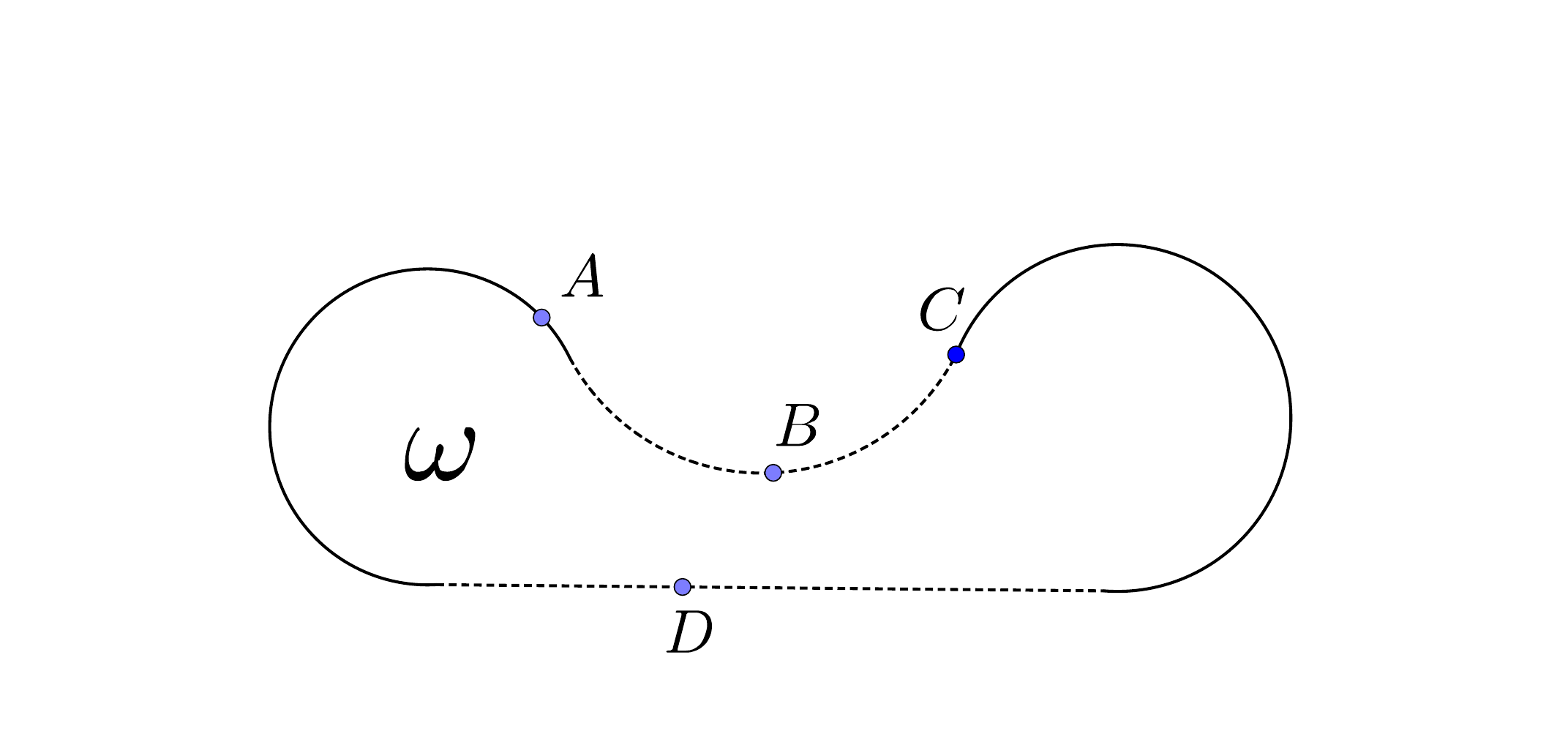}
\caption{Set~$\omega$.}
\label{fig:dom}
\end{figure}

We begin with generalizations of some convex geometry notions. For the classical background, see the book~\cite{Rockafellar}. The material of this subsection is, in some sense, auxiliary. An uninterested reader may skip the proofs without any potential loss of understanding (to make the reading more convenient, we put the proofs of the statements in this subsection into Appendix~\ref{SPG}).
\begin{Def}\label{LocExtrPoint}
Let~$w$ be a subset of~$\mathbb{R}^d$. We call a point~$x \in \partial w$ locally extremal if there are no open segments~$\ell \subset \cl w$ such that~$x \in \ell$. The set of all locally extremal points is called the fixed boundary of~$w$ and denoted by~$\FixedBoundary w$. The set~$\FreeBoundary w = \partial w \setminus \FixedBoundary w$ is called the free boundary.
\end{Def}
For convex sets, locally extremal points are exactly extremal points. For example, for the set~$\Omega$ introduced in Subsection~\ref{s11}, the fixed boundary coincides with~$\partial \Omega_0$ when we impose condition~\eqref{FirstCondition}, see Figure~\ref{fig:strip} also. See Figure~\ref{fig:dom} for visualization; the points~$B, C$ and~$D$ are not locally extremal, whereas the point~$A$ is; the free boundary is marked with a dotted line, whereas the fixed boundary is black.

We have already defined locally concave functions, see Defintion~\ref{LocalConcavity}. 
\begin{Fact}\label{PropertiesOfLCFunctions}
The function~$G$ is locally concave on~$w$ if and only if its restriction to every convex subset of~$w$ is concave.
\end{Fact}

Surely, a pointwise infimum of locally concave functions is locally concave (this follows from the same principle for concave functions). 
\begin{Def}\label{GenLambdaclass}
Let~$w$ be a subset of~$\mathbb{R}^d$\textup, let~$f\colon \FixedBoundary w \to \mathbb{R}$ be a function. Define the class~$\Lambda_{w,f}$ as follows\textup:
\begin{equation*}
\Lambda_{w,f} = \{G\colon w \to \mathbb{R}\cup\{+\infty\}\mid G\hbox{ is locally concave on~$w$}, \quad \forall x \in \FixedBoundary w \quad G(x) \geq f(x)\}.
\end{equation*}
\end{Def}
By the note above, there exists the pointwise minimal function~$\BG = \BG_{w,f}$  defined by formula~\eqref{MinimalLocallyConcave}. Suppose that it nowhere equals~$-\infty$ (for example, this is surely true for the domain~$\Omega$ from Subsection~\ref{s11}), then it belongs to the class~$\Lambda_{w,f}$. Surely,~$\BG_{w,f} = f$ on~$\FixedBoundary w$.  Theorem~$10.1$ in~\cite{Rockafellar} leads us to the following fact.
\begin{Fact}\label{ContinuityInInnerPoints}
Let~$G: w \to \mathbb{R}\cup\{+\infty\}$ be a locally concave function. Then it is continuous at every inner point of~$w$ as a function from~$w$ to~$\mathbb{R}\cup\{+\infty\}$.
\end{Fact}
In particular, if~$\BG$ nowhere equals~$-\infty$, it is continuous. Therefore, the set~$\{x\mid x \in \interior w,\; \BG(x) = \infty\}$ is a relatively closed set; the symbol~$\interior$ denotes the interior. It follows from local concavity that the set~$\{x \mid x \in \interior w,\; \BG(x) = \infty\}$ is relatively open. Therefore, if the set~$\interior w$ is connected,  then~$\BG$ is infinite everywhere on~$\interior w$ if and only if it is infinite in any point.

We need more detailed analysis. We begin with an easy observation.
We will often work with strictly convex sets. We call a convex set~$w$ strictly convex if every point of~$\partial w$ is an exposed point. A point~$x \in \partial w$ is called an exposed point if there exists a hyperplane whose intersection with the closure of~$w$ consists of~$x$ only. 
\begin{Fact}\label{StrictConcavity}
Let~$w$ be a non-empty open convex set. It is strictly convex if and only if~$\FreeBoundary w = \emptyset$\textup, i.e. if~$\partial w$ does not contain segments.
\end{Fact}
\begin{Le}[{\bf Hereditary property}]\label{HereditaryProperty}
Let~$\tilde{w}$ be a closed strictly convex subset of~$w$. Define the function~$\tilde{f}$ on~$\FixedBoundary \tilde{w}$ by the equality
\begin{equation*}
\tilde{f}(x) = \BG(x),\quad x \in \FixedBoundary \tilde{w}. 
\end{equation*} 
In such a case\textup,~$\BG_{w,f}(x) = \BG_{\tilde{w},\tilde{f}}(x)$ for all~$x \in \tilde{w}$.
\end{Le}
\begin{St}\label{ContinuityOnFixedFoundary}
Let~$w$ be a subset of~$\mathbb{R}^d$\textup, let~$x \in \FixedBoundary w$. Suppose that there exists a ball~$B_r(x)$ such that~$B_r(x) \cap w$ is a closed strictly convex set. Suppose~$\BG_{w,f}$ nowhere equals~$+\infty$. Then\textup,~$\BG$ is continuous at~$x$ provided~$f$ is. 
\end{St}
\begin{St}\label{ContinuityOnFreeBoundary}
Let~$x_0$ be a point on the free boundary of~$w$\textup,~$w \subset \mathbb{R}^d$. Suppose there exists some open ball~$B_r(x_0)$ such that~$B_{r}(x_0) \setminus w$ is an open convex set\textup, and~$x_0$ is its exposed point. Suppose that~$\BG < +\infty$ everywhere. Then\textup,~$\BG$ is continuous at the point~$x_0$.
\end{St}
The proofs of Lemma~\ref{HereditaryProperty}, Proposition~\ref{ContinuityOnFixedFoundary}, and Proposition~\ref{ContinuityOnFreeBoundary} are in Appendix~\ref{SPG}.
\subsection{Martingale Bellman function}\label{s22}
We will use minimal amount of probabilistic technique and notation. However, we refer the reader to the book~\cite{Oksendal} for definitions. We are working with discrete-time martingales over an increasing filtration~$\{S_n\}_{n}$ on the standard probability space~$(\mathfrak{S},P)$. For simplicity, all algebras~$S_n$ are finite.
\begin{Def}\label{BasicMartingales}
Let~$w \subset \mathbb{R}^d$ be a closed set. An~$\mathbb{R}^d$-valued martingale~$M$ adapted to~$\{S_n\}_n$ is called an~$w$-martingale if it satisfies the conditions listed below.
\textup{\begin{enumerate}
\item \emph{$S_0 = \{\emptyset, \mathfrak{S}\}$.}
\item \emph{There exists a random variable~$M_{\infty}$ with values in~$\FixedBoundary w$ such that
\begin{equation*}
\E |M_{\infty}| < \infty \quad \hbox{and}\quad M_n = \E(M_{\infty}\mid S_n).
\end{equation*}}
\item \emph{For every~$n \in \mathbb{Z}_+$ and every atom~$\sigma$ in~$S_n$
\begin{equation*}
\conv \Big(M_n(\sigma), \{M_{n+1}(z)\}_{z \in \sigma}\Big) \subset w.
\end{equation*}
The set of all~$w$-martingales is denoted by~$\M_{w}$.}
\end{enumerate}}
\end{Def} 
We note that by L\'evy's zero-one law,~$M_n \to M_{\infty}$ almost surely and in mean.
We give a brief explanation about the third point of the definition above. In particular, it implies that~$M_n$ is a.s. in~$w$. If we had been working with continuous-time martingales, then we could have changed it for the condition ``paths are a.s. inside~$w$''. In such a setting, we consider the set of \^Ito martingales with values in~$w$. They have a.s. continuous paths, therefore, the third point is a consequence of the condition that a martingale is a.s. inside the domain. However, discrete-time martingales do not have continuous paths in any sense, therefore, we need the third condition that forbids a martingale to ``jump over the boundary''. The following lemma shows that~$w$-martingales play the same role for locally concave functions as linear combinations play for concave functions.
\begin{Le}\label{BellmanInduction}
Suppose that~$G$ is locally concave on~$w$ and~$M$ is an~$w$-martingale. Then\textup, the function~$n \mapsto \E G(M_n)$ is non-increasing.
\end{Le}
\begin{proof}
Using the third property of~$w$-martingales for the atom~$\sigma \in S_n$,  we can apply Jensen's inequality to the function~$G$ on the convex set~$\conv \Big(M_n(\sigma), \{M_{n+1}(z)\}_{z \in \sigma}\Big)$ and see that
\begin{equation*}
G\big(M_n(\sigma)\big) \geq \E \big(G (M_{n+1})\,\,\big|\, S_n\big)(\sigma).
\end{equation*}
Averaging, we get
\begin{equation*}
\E G(M_{n+1}) = \E \E\big( G(M_{n+1}) \mid S_n\big) \leq \E G(M_n). 
\end{equation*}
\end{proof}
The procedure just described can be referred to as \emph{the Bellman induction} (for example, see~\cite{Volberg}). To pass to the limit as~$n \to \infty$, we need to consider some subclasses of~$w$-martingales.
\begin{Def}\label{BoundedMartingale}
We say that an~$w$-martingale~$M$ is bounded provided~$M_{\infty}$ is bounded.
\end{Def}
We note that~$|M_n|$ is bounded by the same constant as~$|M_{\infty}|$:
\begin{equation*}
\|M_n\|_{L^{\infty}} = \big\|\E(M_{\infty}\mid S_n)\big\|_{L^{\infty}} \leq \|M_{\infty}\|_{L^{\infty}}.
\end{equation*}
\begin{Def}\label{SimpleMartingale}
We call an~$w$-martingale~$M$ simple if~$M_n = M_{\infty}$ for some~$n$.
\end{Def}
\begin{Def}\label{StronglyMartingaleConnected}
We call~$w \subset \mathbb{R}^d$ a strongly martingale connected domain if for every~$x \in w$ there exists a simple~$w$-martingale starting at~$x$\textup, i.e. a simple~$w$-martingale~$M$ such that~$M_0 = x$.
\end{Def}
\begin{Le}[{\bf Minimal principle}]\label{MaximalPrinciple}
Let~$w$ be a strongly martingale connected domain. Then\textup, for every locally concave function~$G$ on~$w$
\begin{equation*}
\inf_{x\in w}G(x) = \inf_{x \in \FixedBoundary w} G(x).
\end{equation*}
\end{Le}
\begin{proof}
To prove the lemma, it suffices to show that~$G(z) \geq \inf_{x \in \FixedBoundary w} G(x)$ for all~$z \in w$. Let~$M_n^z$ be a simple~$w$-martingale starting at~$z$. Then, by Lemma~\ref{BellmanInduction},~$G(z) \geq \E G(M_n^z)$. Taking~$n$ sufficiently big, we get~$G(z) \geq \E G(M_{\infty}^z) \geq \inf_{x \in \FixedBoundary w} G(x)$, because the values of~$M_{\infty}^z$ are in~$\FixedBoundary w$.
\end{proof}
Now we are ready to define a new Bellman function.
\begin{Def}\label{MartingaleBellman}
Let~$w$ be a strongly martingale connected domain\textup, let~$f$ be a bounded from below function on~$\FixedBoundary w$. Define the martingale Bellman function as
\begin{equation*}
\BM_{w,f}(x) = \sup\big\{\E f(M_{\infty})\,\big|\,\,M_0 = x,\,M \in \M_{w}\big\}, \quad x\in w.
\end{equation*}
If the supremum is taken over the set of bounded~$w$-martingales\textup, then the Bellman function is denoted by~$\BMb$.
\end{Def}
\begin{Rem}\label{BoundedMartingaleNotBoundedFunction}
We note that~$\BMb$ is well defined even for the case where~$f$ is not bounded from below\footnote{We still assume that~$f$ is locally bounded.}.
\end{Rem}
We study this new Bellman function in the next subsection. We refer the reader to Appendix~\ref{SCD}, where we show that~$w$-martingales are worth working on in a broader class of domains than the domains of the type ``convex minus convex''.

\subsection{First duality theorem}\label{s23}
\begin{Le}\label{LocalConcavityOfMartingaleBellman}
For any strongly martingale connected domain~$w$ and any~$f\colon\FixedBoundary w \to \mathbb{R}$ bounded from below we have~$\BM_{w,f} \in \Lambda_{w,f}$.
\end{Le}
\begin{proof}
For every point~$x \in \FixedBoundary w$ there is a constant martingale~$M \equiv x$, therefore,~$\BM \geq f$ on the fixed boundary. Let~$\ell$ be some segment inside~$w$. To verify the concavity of~$\BM\big|_{\ell}$, one has to prove the inequality~$\alpha_+\BM(x_+) + \alpha_-\BM(x_-) \leq \BM(x)$, where~$x = \alpha_+ x_+ + \alpha_- x_-$,~$\alpha_+ + \alpha_- = 1$,~$\alpha_{\pm} \in [0,1]$, and~$[x_+,x_-] \subset \ell$. Fix~$\eta > 0$. Suppose that~$M$ and~$N$ are~$w$-martingales such that
\begin{equation*}
\begin{aligned}
&\BM(x_+) \leq \E f(M_{\infty}) + \eta, \quad &M_0 = x_+;\\
&\BM(x_-) \leq \E f(N_{\infty}) + \eta,\quad &N_0 = x_-. 
\end{aligned}
\end{equation*}
Consider the following martingale~$L$:~$L_0 = x$;~$L_1 = M_0$ with probability~$\alpha_+$,~$L_1 = N_0$ with probability~$\alpha_-$; its parts corresponding to~$M$ and~$N$ develop as~$M$ and~$N$ (i.e.~$P(L_n \in A\mid L_1 = M_0) = P(M_{n-1} \in A)$ and~$P(L_n \in A \mid L_1 = N_0) = P(N_{n-1} \in A)$ for any Borel set~$A$ and any~$n \in \mathbb{N}$). Surely,~$L$ is an~$w$-martingale and~$\E f(L_{\infty}) = \alpha_+\E f(M_{\infty}) + \alpha_-\E f(N_{\infty})$. Thus,
\begin{equation*}
\BM(x) \geq \E f(L_{\infty}) = \alpha_+\E f(M_{\infty}) + \alpha_-\E f(N_{\infty}) \geq \alpha_+ \BM(x_+) + \alpha_-\BM(x_-) - \eta. 
\end{equation*}
Making~$\eta$ arbitrary small, we get the inequality wanted.
\end{proof}
\begin{Rem}\label{LocalConcavityBoundedMartingales}
If~$w$ is a strongly martingale connected domain and~$f$ is locally bounded from below\textup, then~$\BMb$ is locally concave.
\end{Rem}
The Remark is proved by the same argument.
\begin{Le}\label{BellmanInductionLimit}
Let~$w$ be a strongly martingale connected domain\textup, let~$f$ be a bounded from below function on~$\FixedBoundary w$. If~$G \in \Lambda_{w,f}$ is such that~$G$ is continuous at any point of~$\FixedBoundary w$\textup, then\textup,~$G(M_0) \geq \E G(M_{\infty})$ for any~$w$-martingale~$M$.
\end{Le}
\begin{proof}
By Lemma~\ref{BellmanInduction}, the function~$n \mapsto \E G(M_n)$ is non-increasing. We have to prove that~$\underline{\lim} \E G(M_n) \geq \E G(M_{\infty})$. First,~$G(M_n) \to G(M_{\infty})$ a.s., because~$M_n \to M_{\infty}$ a.s. and~$G$ is continuous at the points of the fixed boundary. Second, by Lemma~\ref{MaximalPrinciple},~$G$ is bounded from below. Therefore, the limit relation is a consequence of Fatou's lemma.
\end{proof}
\begin{Rem}
The same assertion is valid if we take~$M$ to be bounded and~$f$ locally bounded from below.
\end{Rem}
\begin{Th}\label{GeometricMartingale}
Let~$w$ be a strongly martingale connected domain\textup, let~$f$ be a bounded from below function on~$\FixedBoundary w$. Suppose that~$\BG_{w,f}$ is continuous at every point of the fixed boundary. Then\textup,~$\BG_{w,f} = \BM_{w,f}$.
\end{Th}
\begin{proof}
The inequality~$\BG \leq \BM$ is contained in Lemma~\ref{LocalConcavityOfMartingaleBellman}, the inequality~$\BM \leq \BG$ follows from Lemma~\ref{BellmanInductionLimit} and Definition~\ref{MartingaleBellman}.
\end{proof}
\begin{Rem}
Conclusion of Theorem~\textup{\ref{GeometricMartingale}} is valid for the  Bellman function~$\BMb_{w,f}$ in place of~$\BM_{w,f}$ even if~$f$ is only locally bounded from below.
\end{Rem}

\section{Functional setting and embedding lemmas}\label{s3}
\subsection{Preliminaries}\label{s31}
Let~$I$ be an interval. Let~$\Omega$ be the same as at the beginning of Subsection~\ref{s11}. The class~$\Classb_{\Omega}$ is the subset of~$\Class_{\Omega}$ consisting of bounded functions. With this class at hand, we can define the Bellman function for  a very general function~$f$,
\begin{equation}\label{BellmanFunctionSimple}
\Bellb_{\Omega,f}(x) = \sup\big\{\av{f(\varphi)}{I}\,\big|\,\,\av{\varphi}{I} = x,\,\varphi \in \Classb_{\Omega}\big\}.
\end{equation}
To define the function~$\Bellb$, we only need~$f$ to be locally bounded from below. The following obvious assertion is a commonplace of the theory: the functions~$\Bell$ and~$\Bellb$ do not depend on the interval~$I$, i.e. if the same classes and functions are constructed on another interval, the resulting Bellman functions are the same (see Remark~$2.1$ of~\cite{IOSVZ} for details). 

We give an alternative form of assumption~\eqref{ThirdCondition}. For every convex set~$\omega \subset \mathbb{R}^2$ and every interior point~$x \in \omega$, there exists the maximal by inclusion convex cone~$C(\omega,x) \subset \omega$ with the vertex~$x$. It is easy to see that this cone is closed. Moreover, it depends on the point~$x$ in a very easy way: if~$x$ and~$y$ are interior points of~$\omega$, then~$C(\omega,y)-y = C(\omega,x)-x$. Therefore, the convex cone~$C(\omega,x) - x$ is independent of the particular choice of~$x$; we call it the maximal inscribed cone of~$\omega$. Assumption~\eqref{ThirdCondition} can be restated using this notation: the set~$\Omega_0$ is infinite and the maximal inscribed cones of~$\Omega_0$ and~$\Omega_1$ are equal.
\begin{Fact}\label{NonemptyClass}
Suppose that~$\Omega$ satisfies assumption~\textup{\eqref{ThirdCondition}}. Then\textup, for every~$x\in \Omega$ there exists a segment~$\ell$ with the endpoints lying on~$\partial \Omega_0$ such that~$x \in \ell \subset \Omega$. 
\end{Fact}
In particular, the fact claims that~$\Omega$ is strongly martingale connected (see Definition~\ref{StronglyMartingaleConnected}). Moreover, it shows that the sets over which the suprema are taken in formulas~\eqref{BellmanFunction} and~\eqref{BellmanFunctionSimple} are non-empty. Indeed, if~$x = \alpha_-x_- + \alpha_+x_+$, where~$\alpha_-+\alpha_+ = 1$,~$\alpha_- \geq 0$,~$\alpha_+ \geq 0$, and~$x_-$ and~$x_+$ are the endpoints of the segment~$\ell$ given by Fact~\ref{NonemptyClass}, then we can define the function~$\varphi_x\colon [0,1]\to\FixedBoundary\Omega$ as follows:
\begin{equation*}
\begin{cases}
\varphi_x(t) = x_-,\quad &t\in[0,\alpha_-];\\
\varphi_x(t) = x_+,\quad &t\in(\alpha_-,1].
\end{cases}
\end{equation*}
It is easy to see that~$\varphi_x\in \Classb_{\Omega}$, because the point~$\av{\varphi_x}{J}$ belongs to the segment~$\ell$ for every~$J$. We note without proof that if assumption~\eqref{ThirdCondition} on~$\Omega$ is violated in the sense that both~$\Omega_0$ and~$\Omega_1$ are infinte, but there exists a ray in~$\Omega_0$ that cannot be translated into~$\Omega_1$, then for some~$x \in \Omega$ the set of functions belonging to~$\Class$ whose average is~$x$, is empty (this is not as obvious as may seem, however, we do not include the proof). Moreover, in such a case~$\Omega$ is not martingale connected. If~$\Omega_0$ is infinte, but~$\Omega_1$ is finite, the situation is interesting again, however, it needs a separate study.

Though the calculation of~$\Bell$ and~$\Bellb$ is a problem, we can compute the values of these functions on the fixed boundary. Indeed, if~$x \in \FixedBoundary \Omega$, then the set of functions we are taking the suprema over consists of one function. Let~$\varphi \in \Class$ be such that~$\av{\varphi}{I} = x$. Let~$\{\Delta_k\}_{k}$ be an arbitrary finite partition of~$I$ into disjoint subintervals. Then,
\begin{equation*}
\av{\varphi}{I} = \sum\limits_{k}\frac{|\Delta_k|}{|I|}\av{\varphi}{\Delta_k}.
\end{equation*}
The convex combination on the right belongs to~$\Omega_0$. By assumption~\eqref{FirstCondition}, if it coincides with~$x$, then all the points~$\av{\varphi}{\Delta_k}$ are equal~$x$. Making~$\{\Delta_k\}_k$ arbitrary small and passing to the limit with the help of the Lebesgue differentiation theorem, we get~$\varphi \equiv x$. Therefore,~$\Bell(x) = f(x)$ and~$\Bellb(x) = f(x)$. This fact can be interpreted as that the Bellman function satisfies the Dirichlet boundary conditions on the fixed boundary.

We end this subsection with an easy auxiliary lemma.
\begin{Le}\label{SummabilityOfModulus}
The function~$\Bell_{|\cdot|}$ is finite. What is more\textup, it is bounded on any compact set~$K$.
\end{Le}
\begin{proof}
We assumed at the very beginning that~$\Omega_0$ does not contain lines. Therefore, there exists some convex cone~$A \subset \mathbb{R}^2$ such that~$\Omega_0 \subset A$. We can find a linear function~$L$\footnote{From now on, we call a linear function what is usually called an affine function, i.e. a function~$f(x)=ax + b$.} such that~$|y| \leq L(y)$ for all~$y \in A$. We can write for any~$\varphi \in \Class$ with the average~$x$:
\begin{equation*}
\av{|\varphi|}{I} \leq \av{L(\varphi)}{I} =L(\av{\varphi}{I}) = L(x).
\end{equation*}
This leads us to the inequality~$\Bell_{|\cdot|}(x) \leq L(x)$, which proves the lemma.
\end{proof}
\subsection{Martingale generates function}\label{s32}
The fixed boundary admits many parametrizations. Fix some line~$\ell$ in the plane such that the orthogonal projection onto~$\ell$ maps~$\FixedBoundary \Omega$ to~$\ell$ bijectively. The inverse map to this projection parametrizes~$\FixedBoundary \Omega$. 
\begin{Def}\label{Monotonicity}
We say that a function~$\varphi\colon I \to \FixedBoundary \Omega$ is monotone if its composition with the projection onto~$\ell$ is monotone. 
\end{Def}
Surely, this notion of monotonicity does not depend on the particular choice of~$\ell$.

For a function~$\varphi\colon I \to\FixedBoundary \Omega$, we denote its distribution, namely, the image of the normalized Lebesgue measure under the mapping~$\varphi$, by~$\mu_{\varphi}$. Thus,~$\mu_{\varphi}$ is a probability measure supported on~$\FixedBoundary \Omega$. We note that for any probability measure~$\mu$ on~$\FixedBoundary \Omega$ there exists a unique mod~$0$ non-increasing function on~$[0,1]$ and a unique mod~$0$ non-decreasing function on~$[0,1]$ whose distribution coincides with~$\mu$. This functions are called the non-increasing and non-decreasing rearrangements of the probability measure~$\mu$ correspondingly.
\begin{Th}\label{MartingaleGeneratesAFucntion}
If~$M \in \M_{\Omega}$ \textup(see Definition~\textup{\ref{BasicMartingales}}\textup)\textup, then the non-decreasing rearrangement of the distribution of~$M_{\infty}$ belongs to~$\Class_{\Omega}$.
\end{Th}
\begin{proof}
Let~$\varphi$ be the non-decreasing rearrangement of the distribution of~$M_{\infty}$. We note that~$\varphi$ is a summable function and~$\av{\varphi}{[0,1]} = M_0$. We have to prove that for any~$J = [a,b] \subset [0,1] = I$ the point~$\av{\varphi}{J}$ is not in~$\Omega_1$. We may assume that~$J \subset \interior I$. Let~$A = \lim_{t \to a_+} \varphi(t)$,~$B = \lim_{t \to b_-}\varphi(t)$. We note that by the monotonicity, for any~$t\in J$ the point~$\varphi(t)$ lies on~$\FixedBoundary \Omega$ between the points~$A$ and~$B$ (we denote the set of all such points by the arc~$\stackrel{\smile}{AB}$). For any function~$f\colon \FixedBoundary \Omega \to \mathbb{R}$ we have
\begin{equation*}
\frac{|J|}{|I|} \av{f(\varphi)}{J}=\E fh_{J}(M_\infty),
\end{equation*}
where
\begin{equation*}
h_{J}=\chi_{\stackrel{\smile}{AB}}-\frac{|\{t\in I \setminus J \mid \varphi(t)=A\}|}{|\{t\in I \mid \varphi(t)=A\}|}\chi_{\{A\}}-\frac{|\{t\in I \setminus J \mid \varphi(t)=B\}|}{|\{t\in I \mid \varphi(t)=B\}|}\chi_{\{B\}}.
\end{equation*}
We note that if the convex hull of~$\stackrel{\smile}{AB}$ belongs to~$\Omega$, then the point~$\av{\varphi}{J}$ is in~$\Omega$ too, because it lies in the mentioned convex hull. So, we assume that the arc~$\stackrel{\smile}{AB}$ is a ``long arc'', i.e. the segment~$AB$ intersects~$\Omega_1$. The idea of the proof is that if the point~$\av{\varphi}{J}$ is not in~$\Omega_1$, then it can be separated from this set by a line. In other words, it suffices to construct a linear function~$L\colon\mathbb{R}^2 \to \mathbb{R}$ such that~$L \geq 0$ on~$\Omega_1$, but~$L(\av{\varphi}{J}) \leq 0$. We consider two cases:~$M_0$ is inside the domain~$S_{A,B}$ and~$M_0$ is outside~$S_{A,B}$ (by~$S_{A,B}$ we mean the closed set shown on Firure~\ref{fig:smile}).

\begin{figure}[h!]
\includegraphics[height=9.5cm]{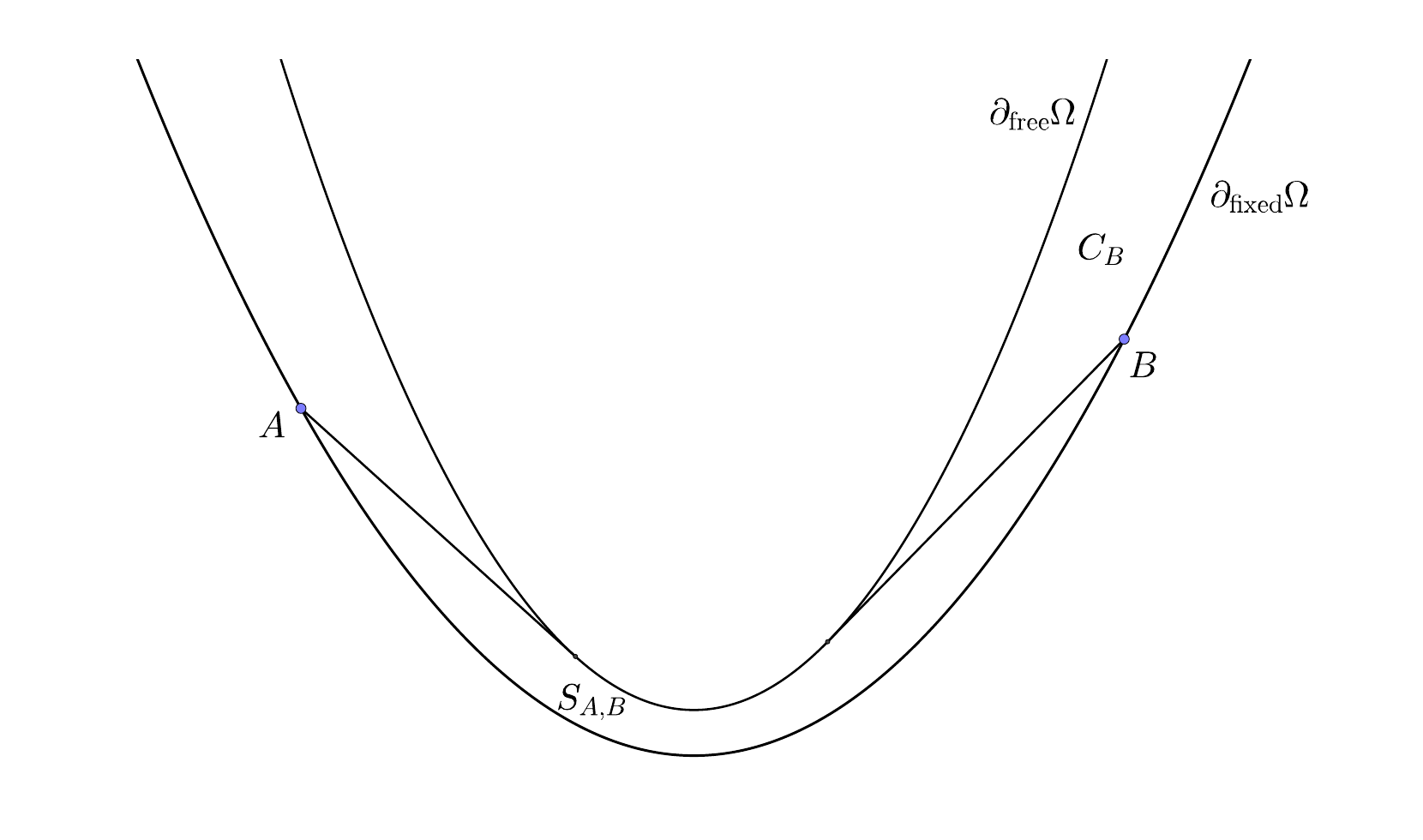}
\caption{Domain~$S_{A,B}$.}
\label{fig:smile}
\end{figure}

\emph{First case.} Let~$L$ be a linear function such that~$L > 0$ on~$\Omega_0 \setminus S_{A,B}$, but~$L(M_0) \leq 0$. For any~$x \in \FixedBoundary \Omega$ we have~$L(x) \geq Lh_J(x)$. Therefore,
\begin{equation*}
\frac{|J|}{|I|}L\big(\av{\varphi}{J}\big)=\frac{|J|}{|I|}\av{L(\varphi)}{J}=\E Lh_J(M_\infty)\leq
\E L(M_\infty)\leq \E L (M_0)\leq 0.
\end{equation*}
The last but one inequality is a consequence of Lemma~\ref{BellmanInductionLimit}.

\emph{Second case.} The domain~$\Omega \setminus S_{A,B}$ consists of two connectivity components. Let~$C_B$ be the one that is adjacent to~$B$. Without loss of generality, we may assume~$M_0 \in C_B$. We take~$L$ such that~$L = 0$ on the common boundary of~$S_{A,B}$ and~$C_B$ and~$L > 0$ on~$\Omega_1$. Surely,~$L \geq Lh_J$ on~$\FixedBoundary \Omega$. Let~$f$ be some continuous function on~$\FixedBoundary \Omega$ such that~$L \geq f \geq Lh_J$ and~$f = 0$ on~$\FixedBoundary \Omega \cap C_B$. Then,
\begin{equation*}
\frac{|J|}{|I|}L\big(\av{\varphi}{J}\big)=\frac{|J|}{|I|}\av{L(\varphi)}{J}=\\
\E Lh_J(M_\infty)\leq \E f(M_\infty)=
\E \BG_{\Omega,f}(M_\infty)\leq \E \BG_{\Omega,f} (M_0) \leq 0.
\end{equation*} 
The last but one inequality is a consequence of Lemma~\ref{BellmanInductionLimit} and Proposition~\ref{ContinuityOnFixedFoundary}. The last inequality follows from the fact that~$\BG_{\Omega,f}$ is non-positive on the common boundary of~$S_{A,B}$ and~$C_B$ (because~$f \leq L$) and zero on~$\FixedBoundary \Omega \cap C_B$, therefore, it is non-positive inside~$C_B$. This implication heuristically follows from the hereditary property, Lemma~\ref{HereditaryProperty}, and the maximal principle~$\sup_{x \in w}\BG_{w,f}(x) = \sup_{x \in \FixedBoundary w} f(x)$. However, the domain whereto we restrict the function is by no means strictly convex, therefore, on the formal level, such an implication does not work. Instead of this, one can assume that~$\BG_{\Omega,f}$ is positive at some point of~$C_B$, then~$\BG_{\Omega,f}$ can be changed for~$\min(\BG_{\Omega,f},0)$ on~$C_B$ and become smaller while retaining to be locally concave (by Fact~\ref{HereditaryOnLines}, see Appendix~\ref{SPG}), which contradicts the minimality. In particular,~$\BG_{\Omega,f} (M_0) \leq 0$.
\end{proof}
\begin{Cor}
If~$M \in \Mb_{\Omega}$\textup, then the monotonic rearrangement of the distribution of~$M_{\infty}$ belongs to~$\Classb$.
\end{Cor}
\subsection{Function generates martingales}\label{s33}
\begin{Def}\label{extension}
Let~$\Omega = \cl \Omega_0 \setminus \Omega_1$ be a domain satisfying  assumptions~\textup{\eqref{FirstCondition}, \eqref{ThirdCondition}}. We call a domain~$\tilde{\Omega}$ an extension of~$\Omega$ if~$\tilde{\Omega} = \cl\Omega_0 \setminus \tilde{\Omega}_1$\textup,~$\tilde{\Omega}_1$ is a convex open set\textup,~$\tilde{\Omega}$ satisfies assumptions~\textup{\eqref{FirstCondition}, \eqref{ThirdCondition}} and~$\cl\tilde{\Omega}_1 \subset \Omega_1$.
\end{Def}
\begin{Th}\label{FunctionGeneratesAMartingale}
Let~$\tilde{\Omega}$ be an extension of~$\Omega$. If~$\varphi \in \Class_{\Omega}$\textup, then there exists a martingale~$M \in \M_{\tilde{\Omega}}$ such that~$M_{\infty}$ is equimeasurable with~$\varphi$.
\end{Th}
\begin{Fact}\label{FunctionDelta}
Define the function~$\Delta\colon \Omega \to [1,\infty)$ as follows\textup:
\begin{equation*}
\Delta(x) = \sup\bigg\{\max\Big(1,\frac{|x-y|}{|x-z|}\Big)\,\bigg|\,\,x \in [y,z],\,y \in \Omega,\,z \in \cl \tilde{\Omega}_1\bigg\}.
\end{equation*}
The function~$\Delta$ is uniformly bounded on every compact subset of~$\Omega$.
\end{Fact}
The following geometric observation goes back to~\cite{Vasyunin} (Lemma~$4$ there) and lies in the heart of the theory. 
\begin{Le}\label{L5}
Let~$\varphi$ be a function on an inteval~$I$ such that~$\varphi \in \Class_{\Omega}$. For any extension~$\tilde{\Omega}$ of~$\Omega$ there exists a partition of~$I$ into two intervals~$I_1$ and~$I_2$ with disjoint interiors such that
\begin{equation*}
[\av{\varphi}{I_1},\av{\varphi}{I_2}] \subset \tilde{\Omega}\quad \hbox{and}\quad \max\Big(\frac{|I_1|}{|I_2|},\frac{|I_2|}{|I_1|}\Big) \leq \Delta(\av{\varphi}{I}).
\end{equation*} 
Here~$\Delta$ is the function introduced in Fact~\textup{\ref{FunctionDelta}}.
\end{Le}
\begin{proof}
For brevity, we assume~$I = [0,1]$ and~$d = \Delta(\av{\varphi}{I})$. Let~$\ell_1(t) = [\av{\varphi}{[0,t]},\av{\varphi}{[0,1]}]$ and~$\ell_2(t) = [\av{\varphi}{[t,1]},\av{\varphi}{[0,1]}]$. For each~$t$, not more than one of the segments~$\ell_1(t)$ and~$\ell_2(t)$ intersects~$\tilde{\Omega}_1$. If~$\ell_2(t)$ intersects~$\tilde{\Omega}_1$, then~$\frac{1-t}{t} < d$. Thus, if~$\frac{1-t_0}{t_0} = d$, then~$\ell_2(t_0)\cap\tilde{\Omega}_1 = \emptyset$. Similarly,~$\ell_1(1-t_0)\cap\tilde{\Omega}_1 = \emptyset$. The set~$A_j$,~$j=1,2$, of all~$t \in [t_0,1-t_0]$ such that~$\ell_j(t)\cap\tilde{\Omega}_1 \ne \emptyset$ is a relatively open subset of~$[t_0,1-t_0]$ that does not cover the whole interval~$[t_0,1-t_0]$. The sets~$A_1$ and~$A_2$ do not intersect, therefore, it follows from the connectivity of~$[t_0,1-t_0]$ that there exists a point~$t$ such that~$t \notin A_1\cup A_2$. We can take~$I_1 = [0,t]$ and~$I_2 = [t,1]$.
\end{proof}
\paragraph{Proof of Theorem~\ref{FunctionGeneratesAMartingale}.} Using Lemma~\ref{L5} inductively, we build a sequence~$\{\{I_k^n\}_{k=1}^{2^n}\}_n$ of partitions of~$I$ such that
\begin{enumerate}
\item For each~$n$ the partition~$\{I_k^{n+1}\}_k$ is a subpartition of~$\{I_k^n\}_k$, moreover, for each~$n$ and~$k$,~ $1 \leq k \leq 2^n$, one has~$I^{n+1}_{2k-1} \cup I^{n+1}_{2k} = I^{n}_k$
\item For each~$n$ and~$k$,~ $1 \leq k \leq 2^{n-1}$, the segment~$\Big[\av{\varphi}{I_{2k-1}^n},\av{\varphi}{I^n_{2k}}\Big]$ lies in~$\tilde{\Omega}$;
\item For each~$n$ and~$k$,~ $1 \leq k \leq 2^n$,~$\max\Big(\frac{|I^{n+1}_{2k-1}|}{|I^{n+1}_{2k}|},\frac{|I^{n+1}_{2k}|}{|I^{n+1}_{2k-1}|}\Big) \leq \Delta(\av{\varphi}{I^n_k})$.
\end{enumerate}
Each partition generates an algebra~$S_n$ of sets. By the first property of our partitions,~$\{S_n\}_n$ is an increasing sequence of algebras. Define the martingale~$M$ by the formula~$M_n = \E(\varphi\big|\,S_n)$. We claim that~$M$ is an $\tilde{\Omega}$-martingale. We have to verify three properties listed in Definition~\ref{BasicMartingales}. The first property is obvious. It follows from the second property of the partitions that~$M$ satisfies the third property of Defintion~\ref{BasicMartingales}. 

To prove the second property, we need to justify that~$\lim_{n \to \infty}\max_{k}|I_k^n| = 0$. Let~$I_n(t)$ denote the unique interval of the~$n$-th partition that contains~$t$ (these functions are defined on the set of full measure). Assume that~$\lim_{n \to \infty}\max_{k}|I_k^n| \ne 0$, then there exists~$t$ such that~$I_n(t) \to J$,~$|J| \ne 0$. In such a case,~$\av{\varphi}{I_n(t)} \to \av{\varphi}{J}$. Therefore, the sequence~$\Delta(\av{\varphi}{I_n(t)})$ is bounded, say, by~$C$, by Fact~\ref{FunctionDelta}. Then, from the third property of the partitions it follows that~$|I_{n+1}(t)| \leq \frac{C}{1+C}|I_n(t)|$, which contradicts the fact that~$I_n(t) \to J$. 

By L\'evy's zero-one law, it follows that~$M_n(t) \to \varphi(t)$ for almost all~$t$. Thus, the second property of Definition~\ref{BasicMartingales} is fulfilled for~$M$ and~$M_{\infty} = \varphi$. \qed
\begin{Rem}
The assertion of Theorem~\textup{\ref{FunctionGeneratesAMartingale}} is not valid\textup, if one asks~$M$ to be in~$\M_{\Omega}$. 
\end{Rem}

\subsection{Corollaries}\label{s34}
\paragraph{Monotonic rearrangements}
\begin{Def}\label{MonotonicRearrangement}
Let~$\varphi\colon I \to \FixedBoundary\Omega$ be a measurable function. A non-decreasing function~$\varphi^*$ \textup(see Defintion~\textup{\ref{Monotonicity}}\textup) that is equimeasurable with~$\varphi$ is called the non-decreasing rearrangement of~$\varphi$.
\end{Def}
\begin{Cor}\label{MonotonicRearrangementCor}
Suppose that~$\Omega$ satisfies conditions~\textup{\eqref{FirstCondition}, \eqref{ThirdCondition}}. If~$\varphi \in \Class_{\Omega}$\textup, then~$\varphi^* \in \Class_{\Omega}$. 
\end{Cor}
\begin{proof}
Assume the contrary, let~$J \subset I$ be an interval such that~$\av{\varphi^*}{J} \notin \Omega$. Let~$\tilde{\Omega}$ be an extension of~$\Omega$ such that~$\av{\varphi^*}{J} \notin \tilde{\Omega}$. We apply Theorem~\ref{FunctionGeneratesAMartingale} to~$\varphi$ and get some~$\tilde{\Omega}$-martingale~$M$ such that~$M_{\infty}$ is equimeasurable with~$\varphi$. The distribution of the non-decreasing rearrangement of~$M_{\infty}$ is equimeasurable with~$\varphi$, therefore, the said monotonic rearrangement coincides with~$\varphi^*$ mod~$0$. Thus, by Theorem~\ref{MartingaleGeneratesAFucntion},~$\varphi^* \in \Class_{\tilde{\Omega}}$, this contradicts the assumption~$\av{\varphi^*}{J} \notin \tilde{\Omega}$.
\end{proof}
\begin{Cor}\label{Concatenation}
Suppose that~$\Omega$ satisfies conditions~\textup{\eqref{FirstCondition}, \eqref{ThirdCondition}}. Let~$\{\varphi_k\}_k$ be a sequence of functions on an interval~$I$\textup, belonging to~$\Class_{\Omega}$. Assume that~$\conv\{\av{\varphi_k}{I}\}_k$ is a bounded subset of~$\Omega$. Let~$\{\alpha_k\}_k$ be a sequence of non-negative numbers such that~$\sum_{k}a_k = 1$. Then\textup, there exists a function~$\varphi \in \Class_{\Omega}$ such that~$\mu_{\varphi} = \sum_k\alpha_k\mu_{\varphi_k}$.
\end{Cor}
\begin{proof}
We take~$\varphi$ to be the non-decreasing rearrangement of the measure~$\sum_k\alpha_k\mu_{\varphi_k}$. We note that~$\av{|\varphi_k|}{I}$ is uniformly bounded by virtue of Lemma~\ref{SummabilityOfModulus}, so~$\varphi \in L^1(I)$. To prove that~$\varphi \in \Class_{\Omega}$, it is sufficient to prove the inclusion~$\varphi \in \Class_{\tilde{\Omega}}$ for all the extensions~$\tilde{\Omega}$ of~$\Omega$. Fix~$\tilde{\Omega}$ and using Theorem~\ref{FunctionGeneratesAMartingale}, construct martingales~$M^k \in \M_{\tilde{\Omega}}$, each on its own probability space~$(\mathfrak{S}_k,P_k)$ with filtration~$\{S^k_n\}_n$, such that~$M^k_{\infty}$ equimeasurable with~$\varphi_k$. Consider a new probability space~$\mathfrak{S} = \bigsqcup_{k=1}^{\infty}\mathfrak{S}_k$,~$P(A) = \sum_{k=1}^{\infty}\alpha_kP_k(A\cap\mathfrak{S}_k)$ for any measurable set~$A$. Define an increasing sequence of algebras~$\{S_n\}_n$ given by formula
\begin{equation*}
A \in S_n \quad \hbox{if and only if} \quad A\cap\mathfrak{S}_k \in S^k_{n-k} \,\, \hbox{for all}\,\, k\leq n \,\,\hbox{and}\,\, A\cap \mathfrak{S}_k \in S^k_0\,\, \hbox{for all}\,\, k > n.
\end{equation*}
Define~$M_n$ by formula
\begin{equation*}
M_n(\sigma) = 
\begin{cases}
M_{n-k}^k(\sigma), \quad & \sigma \in \mathfrak{S}_k, \; k\leq n;\\
\frac{1}{\sum\limits_{m>n}\alpha_m}\sum\limits_{m > n} \alpha_mM^m_0,\quad & \sigma \in \mathfrak{S}_k, \; k> n
\end{cases}
\end{equation*}
with obvious modifications for the case where~$\{\varphi_k\}_k$ is a finite sequence. It is easy to see that~$M$ is a martingale,~$M \in \M_{\tilde{\Omega}}$, and~$M_{\infty}$ is equimeasurable with~$\varphi$. Therefore, by Theorem~\ref{MartingaleGeneratesAFucntion},~$\varphi \in \Class_{\tilde{\Omega}}$.
\end{proof}
\paragraph{Bellman functions}
\begin{Cor}\label{MartingaleLessAnalytic}
Suppose that~$\Omega$ satisfies conditions~\textup{\eqref{FirstCondition}, \eqref{ThirdCondition}}, then
\begin{equation*}
\BMb_{\Omega,f} \leq \Bellb_{\Omega,f}.
\end{equation*}
If~$f$ is bounded from below\textup, then
\begin{equation*}
\BM_{\Omega,f} \leq \Bell_{\Omega,f}.
\end{equation*}
\end{Cor}
This follows from Theorem~\ref{MartingaleGeneratesAFucntion}.
\begin{Cor}\label{AnalyticLessMartingale}
Suppose that~$\Omega$ satisfies conditions~\textup{\eqref{FirstCondition}, \eqref{ThirdCondition}}. For any extension~$\tilde{\Omega}$ of~$\Omega$
\begin{equation*}
\Bellb_{\Omega,f} \leq \BMb_{\tilde{\Omega},f}.
\end{equation*}
If~$f$ is bounded from below\textup, then
\begin{equation*}
\Bell_{\Omega,f} \leq \BM_{\tilde{\Omega},f}.
\end{equation*}
\end{Cor}
This follows from Theorem~\ref{FunctionGeneratesAMartingale}.
\begin{Cor}
The functions~$\Bell$ and~$\Bellb$ are locally concave.
\end{Cor}
This follows from Corollary~\ref{Concatenation} for the case where the sequence~$\{\varphi_k\}_k$ consists of two functions. 
This particular case of Corollary~\ref{Concatenation} for two functions~$\varphi^+$ and~$\varphi^-$ justifies the heuristics that the concatenation of two functions belonging to the class is again in this class (this is not so, but if one takes the monotonic rearrangement of the concatenation, the statement becomes correct). This heuristics was usually used to explain why the searched-for Bellman function is assumed to be locally concave (see, e.g.,~\cite{IOSVZ}, Section~$2.1$).
\paragraph{Distributions of functions belonging to dyadic classes}
\begin{Def}
Let~$Q$ be a cube in~$\mathbb{R}^d$. Let~$\mathcal{D}_Q$ be the set of all dyadic subcubes of~$Q$. Define the dyadic class by the formula
\begin{equation*}
\Classd_{\Omega} = \big\{\varphi\in L^1(Q,\mathbb{R}^2)\,\big|\,\,\varphi(Q) \subset \FixedBoundary\Omega,\ \forall R \in \mathcal{D}_Q \quad \av{\varphi}{R} \in \Omega\big\}.
\end{equation*}
\end{Def}
Dyadic classes ($\BMO^{\mathrm{dyad}}$ and~$A_p^{\mathrm{dyad}}$) are widely studied from different points of view. For the Bellman function approach to similar problems on such classes, see, e.g.~\cite{SV} and~\cite{SV2}.
\begin{Cor}\label{MonotonicDyadic}
There exists an extension~$\tilde{\Omega}$ of~$\Omega$ such that the monotonic rearrangement of any~$\varphi \in \Classd$ belongs to~$\Class_{\tilde{\Omega}}$.
\end{Cor}
We believe that this statement is not new (for particular cases~$\Class = \BMO_{\eps}$ or~$\Class = A_p^{\delta}$; for the case of the Gehring class see~\cite{MN}), however, we have not found it in the literature. For each concrete domain~$\Omega$, the said extension~$\tilde{\Omega}$ can be found by hand (however, the authors do not see a good pattern for its description in full generality). In particular, the procedure below leads to the fact that the monotonic rearrangement of a function belonging to the dyadic~$\BMO$ (or the dyadic Muckenhoupt class) on the cube is in~$\BMO$ (or the Muckenhoupt class).

 We give a sketch of the proof, not going into a detailed study of dyadic classes. Each function from~$\Classd_{\Omega}$ generates a natural martingale in some extension~$\tilde{\Omega}$ of~$\Omega$ (for example, see Lemma~$3.2$ in~\cite{LSSVZ}\footnote{If one follows the constants here, for the particular case of~$\BMO$ we get~$\|f^*\|_{\BMO([0,1])} \leq 2^{\frac{d}{2}} \|f\|_{\BMO^{\mathrm{dyad}}([0,1]^d)}$.}), which, by Theorem~\ref{MartingaleGeneratesAFucntion}, generates a monotonic function in~$\Class_{\tilde{\Omega}}$. It coincides with the monotonic rearrangement of~$\varphi$. Though Corollary~\ref{MonotonicDyadic} may seem very natural, we warn the reader against identifying the dyadic and continuous classes, see Section~$6$ of~\cite{LSSVZ} for some properties that distinguish them. We also mention the forthcoming paper~\cite{SV3} that treats the Bellman function problem on the dyadic classes in arbitrary dimension.

\section{Widening the strip}\label{s4}
\subsection{The plot}\label{s41}
\begin{Th}\label{ExtensionTheorem}
If the domain~$\Omega$ satisfies assumptions~\textup{\eqref{FirstCondition},~\eqref{SecondCondition},~\eqref{ThirdCondition},}~$\FixedBoundary\Omega \in C^2$\textup{,}~$f \in C^2(\FixedBoundary \Omega)$\textup, then
\begin{equation*}
\BG_{\Omega,f} = \inf\big\{\BG_{\tilde{\Omega},f}\,\big|\,\, \tilde{\Omega} \,\,\hbox{is an extension of} \,\,\Omega\big\}.
\end{equation*}
\end{Th}
This theorem finishes the proof of the main theorem for the case of~$C^2$-smooth~$\FixedBoundary\Omega$ and~$f$ (this is discussed in detail in Subsection~\ref{s45}). 
Proof of Theorem~\ref{ExtensionTheorem} occupies the whole present section (except the concluding Subsection~\ref{s45}). For convenience of the reader, we give an informal plot before passing to details. 

If~$\BG_{\Omega,f}$ is infinite, there is nothing to prove. In what follows, we assume the finiteness of this function to avoid confusion. Our aim is, with any positive~$\delta$ at hand, to construct an extension~$\tilde{\Omega}_{\delta}$ and a locally concave function~$\tilde{G}_{\delta}$ on it such that
\begin{equation*}
\forall x \in \Omega \quad\BG_{\Omega,f}(x) \leq \tilde{G}_{\delta}(x) \leq \BG_{\Omega,f}(x) + \delta.
\end{equation*} 
This is done in several steps. We begin with an investigation on the behavior of a minimal locally concave function near the free boundary. By Proposition~\ref{ContinuityOnFreeBoundary},~$\BG$ is continuous at the points of~$\FreeBoundary \Omega$. If the boundary and the boundary conditions are~$C^1$-smooth, then at any point~$x \in \FreeBoundary \Omega$ there exists a linear function~$L[\BG_{\Omega,f};x]$ such that~$L[\BG_{\Omega,f};x](x) = \BG(x)$ and~$\BG(y) \leq L[\BG_{\Omega,f};x](y)$ if the points~$x$ and~$y$ see each other\footnote{By this we mean that~$[x,y] \subset \Omega$.}. In other words, the graph of~$L[\BG_{\Omega,f};x]$ is a supporting plane at the point~$(x,\BG(x))$ to the graph of~$\BG$ restricted to the set of points that see~$x$.

Though the graph of~$L[\BG_{\Omega,f};x]$ is a supporting plane, the value of~$\BG$ can be bigger than the value of~$L[\BG_{\Omega,f};x]$ at the point~$y$ that is near~$x$ but is not seen from it (e.g.,~$y$ is another point of~$\FreeBoundary \Omega$), and this is very common for minimal locally concave functions. However,~$L[\BG_{\Omega,f};x]$ can be chosen in such a way that the gain the value~$\BG(y)$ can obtain ``against the usual concavity'' is small,~$\BG(y) - L[\BG_{\Omega,f};x](y) \leq O(|x-y|^3)$. From this inequality, it is not hard to see that we can find some small locally concave on~$\Omega$ function~$E_{\delta}$ such that~$0 \leq E_{\delta} \leq \frac{\delta}{2}$ and the sum~$\BG + E_{\delta}$ is strictly locally concave. This material constitutes Subsection~\ref{s42}.

We have to construct some locally concave function on some extended domain. The idea is to extend the function~$\BG + E_{\delta}$ through the free boundary (i.e. construct the function~$\tilde{G}_{\delta}$ on~$\tilde{\Omega}_{\delta}$ in such a way that~$\tilde{G}_{\delta} = \BG + E_{\delta}$ on~$\Omega$). In Subsection~\ref{s44}, we prove that any strictly locally concave function can be extended in such a way, thus proving Theorem~\ref{ExtensionTheorem}. The procedure may seem a bit illogical at the first sight: minimal locally concave functions usually cannot be extended through the free boundary. However, they can be perturbed to become strictly locally concave, and only then extended.

Before turning to the proof, we mention an easy consequence of the extension theorem.
\begin{Cor}\label{ExtensionForSmoothFunctions}
The statement of Theorem~\textup{\ref{ExtensionTheorem}} is valid with milder condition~$f \in C(\FixedBoundary \Omega)$ instead of~$f \in C^2(\FixedBoundary\Omega)$. 
\end{Cor}
\begin{proof}
For any~$\eps > 0$ we can find a function~$f_{\eps}\in C^2(\FixedBoundary \Omega)$ such that~$0 \leq f_{\eps} - f < \eps$. This leads to~$0 \leq \BG_{f_{\eps}} - \BG_{f} \leq \eps$. Applying Theorem~\ref{ExtensionTheorem} to the function~$f_{\eps}$ and then sending~$\eps$ to zero, we prove the corollary.
\end{proof}

\subsection{Superdifferential and local growth along free boundary}\label{s42}
\begin{figure}[h!]
\includegraphics[height=9.5cm]{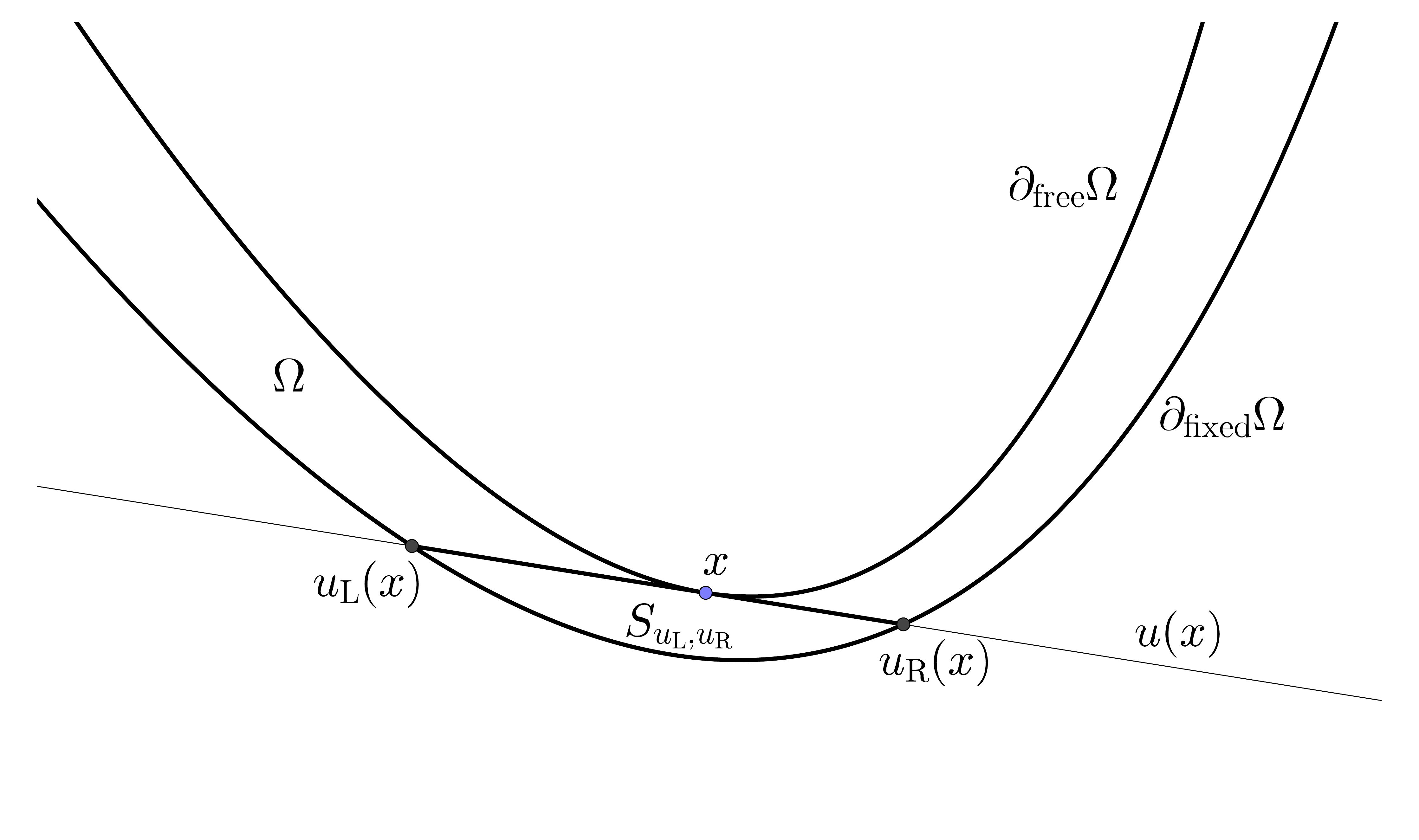}
\caption{Illustration to the notation.}
\label{fig:ill}
\end{figure}

We introduce some notation. Let~$x$ be a point in~$\FreeBoundary \Omega$. Define~$\ul(x)$ and~$\ur(x)$ to be two points on~$\FixedBoundary \Omega$ such that~$x \in [\ul,\ur]$ and this segment touches~$\FreeBoundary \Omega$. Such points exist due to assumptions~\eqref{SecondCondition} and~\eqref{ThirdCondition}. The line that passes through~$\ul$ and~$\ur$ is denoted by~$u = u(x)$. The closure of the part of~$\Omega$ that is separated from~$\Omega_1$ by the line~$u$ is called~$S_{\ul,\ur}$.  We note that~$S_{\ul,\ur}$ consists of those points in~$\Omega$ that see~$x$.  We parametrize the fixed boundary with the inverse of some orthogonal projection as we did at the beginning of Subsection~\ref{s32}. With this projection at hand, we say that one point lies on the left of another if it has smaller projection (we also fix some orientation of the line whereto we project). We assume that for any~$x$ the point~$\ul(x)$ lies on the left of~$\ur(x)$. See Figure~\ref{fig:ill} for clarity.

\begin{Le}\label{Superdifferential}
Suppose that~$\Omega$ and~$f$ satisfy the assumptions of Theorem~\textup{\ref{ExtensionTheorem}} and~$x \in \FreeBoundary \Omega$. There exists a linear function~$L$ such that~$L(x) = \BG(x)$ and~$L(y) \geq \BG(y)$ whenever~$[x,y] \subset \Omega$.
\end{Le}
\begin{proof}
The function~$\BG\big|_{[\ul(x),\ur(x)]}$ is concave, therefore, there is some linear function~$\ell\colon u \to \mathbb{R}$ such that~$\ell(x) = \BG(x)$ and~$\ell(y) \geq \BG(y)$ for any~$y \in [\ul,\ur]$. The function~$L$ we are searching for will coincide with~$\ell$ on~$[\ul,\ur]$. 

By differentiability of~$f$ on the fixed boundary, there exist positive~$\eps$ and~$c$ such that
\begin{equation*}
|f(z) - f(\ul)| \leq c\dist(z,u)\quad \hbox{for}\,\,\, z \in \FixedBoundary \Omega \cap S_{\ul,\ur}, \,\dist(z,\ul) < \eps,
\end{equation*}
and a similar inequality holds true with~$\ur$ instead of~$\ul$. Therefore, there exists some linear function~$L_1$ that coincides with~$\ell$ on~$u$ and satisfies the inequality~$f(z) \leq L_1(z)$ for all~$z$ such that~$z \in\FixedBoundary \Omega \cap S_{\ul,\ur}$ and~$\dist(z,u) < \eps$.

On the other hand, the function~$f$ is bounded on the set~$S_{\ul,\ur} \cap\FixedBoundary\Omega \cap \{z \mid \dist(z,u) \geq \eps\}$, thus there exists a linear function~$L_2$ such that it coincides with~$\ell$ on~$u$ and~$L_2(z) \geq f(z)$ for~$z \in \FixedBoundary \Omega \cap S_{\ul,\ur}$ such that~$\dist(z,u) \geq \eps$.

Let~$L$ be the maximum of~$L_1$ and~$L_2$ (on~$S_{\ul,\ur}$). We see that~$L(z) \geq f(z)$ on~$\FixedBoundary \Omega \cap S_{\ul,\ur}$ and on~$[\ul,\ur]$. By the minimality of~$\BG$,~$L(y) \geq \BG(y)$ for all~$y \in S_{\ul,\ur}$. Indeed, if this inequality does not hold, then the function that equals~$\min(\BG,L)$ on~$S_{\ul,\ur}$ and~$\BG$ everywhere else is smaller than~$\BG$, satisfies the boundary conditions, and is locally concave (all this contradicts the minimality of~$\BG$).  
\end{proof}
We have built a function~$L$ that coincides with some linear function~$\ell$ on~$u$ and exceeds~$\BG$ on~$S_{\ul,\ur}$. Surely, there exists a pointwise minimal on~$S_{\ul,\ur}$ function among those which satisfy the same conditions (with the function~$\ell$ fixed). We call it~$L[\BG_{\Omega,f};x]$. Note that its choice depends on our choice of~$\ell$\footnote{One can prove that~$\BG\big|_{(\ul,\ur)}$ is in~$C^1$ with the same regularity assumptions at hand. Therefore, the function~$\ell$ is unique. We do not need this.}. Let us clarify the notation. If~$\omega$ is a set,~$G$ is a locally concave function on it, then~$L[G;x]$ denotes a linear function with two properties: first,~$L[G;x](x) = G(x)$; second,~$L[G;x](y) \geq G(y)$ for any~$y$ in some neighborhood of~$x$,~$[x,y] \subset \omega$. The importance of collections of such linear functions is emphasized by Fact~\ref{SuperDifferentialConcavity} below. We note that the function~$L[G;x]$ may not be well defined. However, we will always clarify what specific linear function we choose.

It follows from the minimality of~$L[\BG;x]$ that for each~$x$ there exists some point~$e_x \in S_{\ul(x),\ur(x)} \cap \FixedBoundary \Omega$ such that~$L[\BG;x](e_x) = f(e_x)$. The following fact is an easy consequence of the~$C^2$-smoothness assumptions.
\begin{Fact}\label{LocalMaxima}
In the assumptions of Theorem\textup{~\ref{ExtensionTheorem},} for any~$x\in\FreeBoundary \Omega$\textup, there exists a point~$e_x$ inside the set~$S_{\ul(x),\ur(x)} \cap \FixedBoundary \Omega$ such that for any~$z$ in~$\FixedBoundary \Omega$\textup,~$[x,z] \subset \Omega$\textup, the following equality holds\textup:
\begin{equation*}
|L[\BG;x](z) - f(z)| = O(|z - e_x|^2),
\end{equation*}
where the~$O$ is uniform when~$x$ is in some compact set.
\end{Fact} 
It follows from the construction of the linear functions~$L[\BG;x]$ that they are uniformly bounded on compact sets (uniformly with respect to~$x$ lying in some compact set). Therefore, the slopes of these linear functions are uniformly bounded.
\begin{St}\label{SuperConcavityEstimate}
For each compact set~$K \subset \FreeBoundary \Omega$\textup, there exist constants~$\eps > 0$ and~$c > 0$ such that
\begin{equation*}
\BG(y) \leq L[\BG;x](y) + c|x-y|^3 \quad {whenever}\ x\in K, \,\,y\in \FreeBoundary\Omega,\,\,|x-y| < \eps.
\end{equation*} 
\end{St}
\begin{proof}
Assume that~$y$ lies on the right of~$x$. If~$\eps$ is sufficiently small, then~$\ul(y) \in S_{\ul(x),\ur(x)}$ and~$\ur(x) \in S_{\ul(y),\ur(y)}$. We consider two cases:~$[e_x,y] \subset \Omega$ and~$[e_x,y] \not\subset \Omega$. Here~$e_x$ is the point defined in Fact~\ref{LocalMaxima}.

If the segment~$[e_x,y]$ lies inside~$\Omega$, then~$e_x \in S_{\ul(y),\ur(y)}$. Consider the function~$\BG - L[\BG;x]$ restricted to the segment~$[e_x,y]$. This function is concave, it has zero value at the point~$z = e_x$ and non-positive value at the point~$[\ul(x),\ur(x)] \cap [e_x,y]$. Therefore, it has non-positive value at the point~$y$, and the inequality~$\BG(y) \leq L[\BG;x](y)$ holds without the additional term on the right-hand side\footnote{In the case where~$e_x = \ur(x)$ this does not work exactly. In such a case, one has to approximate~$e_x$ by some point~$z \in \FixedBoundary\Omega \cap S_{\ul(x),\ur(x)}$, do the same reasoning for the segment~$[y,z]$, use Fact~\ref{LocalMaxima}, and passing to the limit as~$z \to \ur(x)$ get the same inequality without the additional term.}.

If~$e_x \notin S_{\ul(y),\ur(y)}$, then we use Fact~\ref{LocalMaxima} and see that
\begin{equation*}
\big|f(\ul(y)) - L[\BG;x](\ul(y))\big| = O(|\ul(y) - e_x|^2) = O(|\ul(x) - \ul(y)|^2) = O(|x-y|^2).
\end{equation*}
Now we consider the function~$\BG - L[\BG;x]$ restricted to the segment~$[\ul(y),y]$. This function is bounded from below by~$O(|x-y|^2)$ at the point~$\ul(y)$ and it is non-positive at the point~$p_{x,y} = [\ul(y),y] \cap [\ul(x),\ur(x)]$. Therefore, at the point~$y$, it is bounded from above by the value
\begin{equation*}
\frac{|p_{x,y} - y|}{|p_{x,y} - \ul(y)|}O(|x-y|^2) = O(|x-y|^3).
\end{equation*}
Now one only has to look through the proof again and convince himself that the constants are uniform on the compact set~$K$.
\end{proof}
By a strictly concave~$C^2$-smooth function we mean a function whose Hessian (the matrix of second differential) is strictly positive-definite at every point. 
\begin{Le}\label{SmallStrictlyConcaveFunction}
For any convex domain~$\Omega_0 \subset \mathbb{R}^2$ that does not contain lines and any~$\delta > 0$ there exists a strictly concave~$C^2$-smooth function~$E_{\delta}$ on~$\cl\Omega_0$ such that~$0 < E_{\delta}(x) < \delta$ for any~$x \in \cl\Omega_0$.
\end{Le}
\begin{proof}
If~$\Omega_0$ does not contain lines, then it lies inside some convex cone. Therefore, there exist non-parallel lines~$\ell_1$ and~$\ell_2$ such that~$\ell_1\cap\cl\Omega_0 = \ell_2\cap\cl\Omega_0 = \emptyset$. One can easily verify that the function
\begin{equation*}
E_{\delta}(x) = \frac{\delta}{2}\big(2 - e^{-\dist(x,\ell_1)} - e^{-\dist(x,\ell_2)}\big)
\end{equation*}
does not exceed~$\delta$, is positive, and is strictly concave on~$\cl\Omega_0$.
\end{proof}
\begin{Cor}\label{StrictlyConcaveUpperApproximation}
Suppose that~$\Omega$ and~$f$ satisfy the assumptions of Theorem~\textup{\ref{ExtensionTheorem}}\textup; let~$E_{\delta}$ be the function constructed in the previous lemma. For any~$\delta > 0$ the function~$\BG_{\Omega,f} + E_{\delta}$ has the properties listed below.
\textup{\begin{enumerate}
\item~$\BG_{\Omega,f} \leq \BG_{\Omega,f} + E_{\delta} \leq \BG_{\Omega,f} + \delta$.
\item \emph{For any point~$x \in \FreeBoundary \Omega$ there exists a linear function~$L[\BG_{\Omega,f} + E_{\delta};x]$ such that~$\BG_{\Omega,f}(y) + E_{\delta}(y) \leq L[\BG_{\Omega,f} + E_{\delta};x](y)$ for all~$y \in S_{u_1(x),u_2(x)}$\textup, but~$\BG_{\Omega,f}(x) + E_{\delta}(x) = L[\BG_{\Omega,f} + E_{\delta};x](x)$\textup; for any compact set~$K \subset \FreeBoundary \Omega$ the value~$\sup_{x\in K} \big(|\nabla L[\BG_{\Omega,f} + E_{\delta};x]| + |L[\BG_{\Omega,f} + E_{\delta};x](0)|\big)$ is finite.}
\item \emph{For any compact set~$K\subset\FreeBoundary \Omega$ there exist~$\eps_K > 0$ and~$c_K > 0$ such that
\begin{equation*}
\BG_{\Omega,f}(y) + E_{\delta}(y) \leq L[\BG_{\Omega,f} + E_{\delta};x](y) -c_K|x-y|^2\,\,\hbox{whenever}\,\, x \in K,\,\, y\in \FreeBoundary\Omega, \,\,|x-y| < \eps_K.
\end{equation*}} 
\end{enumerate}}
\end{Cor}
\begin{proof}
The first condition follows from the construction, the second one is a consequence of Lemma~\ref{Superdifferential} (the functions~$L[\BG_{\Omega,f} + E_{\delta};x]$ are defined as the~$L[\BG_{\Omega,f};x]$ plus the linear functions that approximate~$E_{\delta}$ at the same points), and the third one is provided by Proposition~\ref{SuperConcavityEstimate} and the strict concavity of~$E_{\delta}$. 
\end{proof}
\subsection{Extension of strictly locally concave function}\label{s44}
We are going to extend the function~$\BG_{\Omega,f} + E_{\delta}$. What we need to construct an extension are the linear functions~$L[\BG_{\Omega,f} + E_{\delta};x]$ only. We treat~$\FreeBoundary \Omega$ as a fixed boundary and try to construct a domain~$\hat{\Omega}_{\delta}$,~$\FreeBoundary\Omega = \FixedBoundary \hat{\Omega}_{\delta}$, and a locally concave function~$\hat{G}_{\delta}$ on it such that~$\hat{G}_{\delta} = \BG_{\Omega,f} + E_{\delta}$ on~$\FixedBoundary\hat{\Omega}$ and for each~$z$ the inequality~$\hat{G}_{\delta}(z) \leq L[\BG_{\Omega,f} + E_{\delta};x](z)$ holds, whenever~$z$ and~$x$ see each other in~$\hat{\Omega}_{\delta}$ and~$x \in \FixedBoundary \hat{\Omega}_{\delta}$. Then we can glue the functions~$\hat{G}_{\delta}$ and~$\BG_{\Omega,f} + E_{\delta}$ to get~$\tilde{G}_{\delta}$ (a locally concave extension of~$\BG_{\Omega,f} + E_{\delta}$) defined on~$\tilde{\Omega}_{\delta} = \Omega \cup \hat{\Omega}_{\delta}$. The local concavity of~$\tilde{G}_{\delta}$ will follow immediately from a simple fact indicated below (which is a consequence of a similar principle for concave functions on an interval).
\begin{Fact}\label{SuperDifferentialConcavity}
Let~$w$ be a subset of~$\mathbb{R}^d$\textup, let~$G$ be a function on it. If for each point~$x \in w$ there exists a linear function~$L[G;x]$ such that~$L[G;x](x) = G(x)$ and~$L[G;x](y) \geq G(y)$ for all~$y$ in some neighborhood of~$x$\textup,~$[x,y] \subset w$\textup, then~$G$ is locally concave on~$w$.
\end{Fact}

We need an abstract lemma about extending domains.
\begin{Le}\label{FirstExtension}
Suppose that~$\Omega_0$ is a strictly convex unbounded open non-empty subset of~$\mathbb{R}^2$ that does not contain lines. For any~$\rho\colon \partial \Omega_0 \to (0,1]$ that is separated from zero on each compact subset of~$\partial \Omega_0$\textup, there exists a convex open set~$\Omega_1$ such that~$\Omega = \cl \Omega_0 \setminus \Omega_1$ satisfies assumptions~\textup{\eqref{FirstCondition}} and {\eqref{ThirdCondition}} and whenever the points~$x$ and~$y$ belonging to $\partial \Omega_0$ see each other in~$\Omega$\textup, the inequality~$|x-y| \leq \rho(x)$ holds true.
\end{Le}
\begin{proof}
Define~$\Omega'_1$ by the formula
\begin{equation*}
\Omega'_1 = \conv\bigg(\Big\{\frac{x + y}{2}\,\bigg|\,\, x,y \in \partial\Omega_0,\,|x-y| > \frac{\rho(x)}{2}\Big\}\bigg).
\end{equation*}  
It follows from the strict convexity of~$\Omega_0$ that the closure of~$\Omega'_1$ does not intersect~$\partial \Omega_0$. 
We can enlarge the set~$\Omega'_1$ up to some bigger strictly convex open set~$\Omega_1$ still lying strictly inside~$\Omega_0$  such that~$\Omega = \cl \Omega_0 \setminus \Omega_1$ satisfies condition~\eqref{FirstCondition} and~\eqref{ThirdCondition}. 
We only have to prove that the points~$x$ and~$y$ from~$\partial \Omega_0$ do not see each other in~$\Omega$ if~$|x-y| > \rho(x)$. This follows from the fact that the midpoint of~$[x,y]$ is in~$\Omega'_1$ provided~$|x-y| > \rho(x)$. 
\end{proof}
\begin{Cor}\label{StrictlyConcaveMajorant}
Suppose that~$\Omega_0$ is a strictly convex unbounded open non-empty subset of~$\mathbb{R}^2$\textup, which has~$C^2$-smooth boundary and does not contain lines. Suppose the linear functions~$L_x$\textup,~$x\in\partial\Omega_0$\textup, satisfy the condition that for any compact set~$K \subset \partial\Omega_0$ there exist positive numbers~$\eps_K$ and~$c_K$ such that for any~$x 
\in K$ and~$y\in \partial\Omega_0$\textup,
\begin{equation}\label{LightBound}
L_y(y) \leq L_x(y) - c_K|x-y|^2
\end{equation}
whenever~$|x-y| < \eps_K$.
In such a case\textup, there exists a domain~$\Omega_1$ such that~$\Omega = \cl\Omega_0 \setminus \Omega_1$ satisfies conditions~\eqref{FirstCondition} and~\eqref{ThirdCondition} and inequality~\textup{\eqref{LightBound}} holds when~$x \in K$ and~$y\in \partial\Omega_0$ see each other in~$\Omega$ \textup(with some different parameter~$c_K'$\textup)\textup; the domain~$\Omega$ also satisfies the condition that the points see each other in it only when the distance between them does not exceed~$1$.
\end{Cor}
\begin{proof}
This corollary is a particular case of Lemma~\ref{FirstExtension}. Consider a locally finite covering of~$\partial \Omega_0$ by a family~$\mathfrak{K}$ of compact sets. For each~$x \in \partial \Omega_0$, we fix a compact set~$\varkappa_x \in \mathfrak{K}$ such that~$x \in \varkappa_x$. Define the function~$\rho$ by the formula
\begin{equation*}
\rho(x) = \min\Big(0.1, \frac{\eps_{\varkappa_x}}{2}\Big),\quad x \in \partial \Omega_0.
\end{equation*} 
This function is surely separated from zero on any compact set (because it attains only a finite number of values on this set). We apply Lemma~\ref{FirstExtension} with this particular function~$\rho$ to build the domain~$\Omega$. We have to verify inequality~\eqref{LightBound} for any compact set~$K \subset \partial \Omega_0$. It follows from the very definition of the objects that this inequality holds with the parameters~$C_K'$ chosen by the rule
\begin{equation*}
c'_K = \min\Big\{c_{\varkappa}\,\Big|\,\, \varkappa \in \mathfrak{K}, \,\, \varkappa \cap K \ne \emptyset\Big\}.
\end{equation*}
\end{proof}
Let~$\Omega_0$ and the collection of linear functions~$\{L_x\}_x$ be as in the previous corollary. We construct the function~$G$ by the formula
\begin{equation}\label{LightFunction}
G(z) = \inf\{L_x(z)\big|\,\,[z,x] \subset \Omega,\,\,x \in \FixedBoundary\Omega\},
\end{equation}
where the domain~$\Omega$ is the one constructed in the previous corollary. This function, surely, satisfies~$G(x) = L_x(x)$  and~$G(z) \leq L_x(z)$ whenever~$x \in \FixedBoundary\Omega$ and~$[z,x] \subset \Omega$. We are going to prove its local concavity on some domain~$\Omega^*$ which is a subdomain of~$\Omega$. We almost obtain it from formula~\eqref{LightFunction} and Fact~\ref{SuperDifferentialConcavity}: if~$G(z) = L_{\check{z}}(z)$ and~$\check{z}$ sees some neighborhood of~$z$, then~$G$ is smaller than~$L_{\check{z}}$ in some neighborhood of~$z$  by the definition of~$G$. It occurs that if~$z$ is not far from~$\FixedBoundary\Omega$, then~$\check{z}$ is also not far from~$z$, and thus sees some neighborhood of~$z$. Here is the precise statement.
\begin{Le}\label{ConcavityNearTheBoundary}
Suppose that~$\Omega$ satisfies the assumptions~\textup{\eqref{FirstCondition}, \eqref{ThirdCondition},} and that if~$[x,y] \subset \Omega$\textup, then~$\hbox{$|x-y| \leq 1$}$. Suppose the linear functions~$L_x$\textup,~$x \in \FixedBoundary\Omega$\textup, satisfy the condition that for any compact set~$K \subset \FixedBoundary \Omega$ there exists~$c_K > 0$ such that for any points~$x\in K$\textup,~$y\in\FixedBoundary\Omega$ that see each other in~$\Omega$ inequality~\textup{\eqref{LightBound}} holds true. We also assume that~$|\nabla L_x|+ |L_{x}(0)|$ is uniformly bounded when~$x$ belongs to a compact set.
If the function~$G$ is given by formula~\textup{\eqref{LightFunction},} then for any compact set~$K \subset \FixedBoundary\Omega$ there exists~$\eps_K > 0$ such that if~$x \in K$\textup,~$y\in\FixedBoundary\Omega$\textup, and~$|x-y| \leq \eps_K$\textup, then~$[x,y] \subset \Omega$ and for any~$z \in [x,y]$ there exists a linear function~$L[G;z]$ such that~$G(z) = L[G;z](z)$\textup, but~$G$ is not bigger than~$L[G;z]$ in some neighborhood of~$z$.
\end{Le}
\begin{proof}
We need to fix two parameters for each compact set~$K \subset \FixedBoundary \Omega$,~$d_K = \dist(K,\FreeBoundary \Omega)$ and~$s_K = \sup|\nabla L_x|$, where the supremum is taken over all~$x \in \FixedBoundary\Omega$ that are seen from some point of~$K$ (this set is compact due to our assumptions on~$\Omega$). The parameter~$\eps_K$ will satisfy several inequalities to be specified later, the first one is that~$\eps_K$ is much smaller than~$d_k$ (say,~$\eps_K \leq 0.1d_K$).

We are going to prove that there exists~$\epsilon_K > 0$ such that for any points~$z_0, p \in K$ such that~$z_0$ and~$p$ see each other in~$\Omega$ and any point~$z \in \Omega$ such that~$|z-z_0| \leq \epsilon_K$,~$|z-p| \geq 0.5d_K$, the inequality~$L_{z_0}(z) \leq L_p(z)$ holds. This is easy:
\begin{equation*}
\begin{aligned}
L_{z_0}(z) \leq L_{z_0}(z_0) + s_K|z-z_0| \stackrel{\hbox{\tiny \eqref{LightBound}}}{\leq} L_p(z_0) - c_K|z_0 - p|^2 + s_K|z-z_0| \leq\\ L_p(z) -c_K|z_0-p|^2 + 2s_K|z-z_0| \leq
L_p(z) -c_K\big(|z-p| - |z_0-z|\big)^2 + 2s_K|z-z_0| \leq\\ L_p(z) - c_K(0.5d_K - \epsilon_K)^2 + 2s_K\epsilon_K \stackrel{\scriptscriptstyle\epsilon_K \leq 0.2d_K}{\leq} L_p(z) - 0.09c_Kd_K^2 + 2s_K\epsilon_K \stackrel{\scriptscriptstyle\epsilon_K \leq 0.04\frac{c_Kd_K^2}{s_K}}{\leq} L_p(z).
\end{aligned}
\end{equation*}

So, for any compact set~$K$ we take~$\epsilon_K = \min(0.1d_K, 0.04\frac{c_Kd_K^2}{s_K})$. Now we take~$\eps_K = 0.1\epsilon_{K_1}$, where~$K_1$ is the set of all points on~$\FixedBoundary \Omega$ that are seen from the points that are seen from~$K$. In such a case,
\begin{equation*}
G(z) = \inf\{L_p(z)\,|\,\,|z - p| \leq 0.5d_{K_1}, \,\, p \in \FixedBoundary \Omega\}
\end{equation*}
if~$z \in [x,y]$,~$|x-y| \leq \eps_K$,~$x \in K$,~$y \in \FixedBoundary \Omega$. Indeed, choose any point~$p \in \FixedBoundary \Omega$ such that~$|z-p| \geq 0.5 d_{K_1}$ and~$p$ sees~$z$. In such a case,~$p$ sees either $x$ or $y$. We take $z_0=x$ or $z_0=y$ in such a way that~$p$ sees~$z_0$. Then, both~$p$ and~$z_0$ are in~$K_1$,~$|z-z_0|\leq \eps_K \leq \epsilon_{K_1}$, and~$|z-p|\geq 0.5d_{K_1}$. Thus, by the argument above,~$L_p(z)\geq L_{z_0}(z) \geq \min(L_x(z),L_y(z))$.

The set of points~$p$ that see~$z$ and satisfy~$|z-p| \leq 0.5d_{K_1}$ is compact, thus there exists~$\check{z}$ such that~$G(z) = \lim_{z_n \to \check{z}}L_{z_n}(z)$. Take~$L[G;z]$ to be~$\lim_{z_n \to \check{z}}L_{z_n}$ (the sequence~$\{L_{z_n}\}_n$ is uniformly bounded as a set of linear functions, therefore, it has some convergent subsequence; without loss of generality, we assume that~\{$L_{z_n}\}_n$ has a limit itself). Surely,~$\check{z}$ sees some neighborhood of~$z$ (and so does~$z_n$ eventually), because~$|z-\check{z}| \leq 0.5 d_{K_1}$ and~$\check{z} \in K_1$. Therefore, by the construction of~$G$ (formula~\eqref{LightFunction}),~$G$ does not exceed~$L[G;z]$ on this neighborhood.
\end{proof}
\begin{St}\label{LightExtension}
Suppose that~$\Omega$ satisfies the assumptions~\textup{\eqref{FirstCondition}} and~\textup{\eqref{ThirdCondition},} and that if~$[x,y] \subset \Omega$\textup, then~$\hbox{$|x-y| \leq 1$}$. Suppose the linear functions~$L_x$\textup,~$x \in \FixedBoundary\Omega$\textup, satisfy the condition that for any compact set~$K \subset \FixedBoundary \Omega$ there exists~$c_K > 0$ such that for any points~$x\in K$\textup,~$y\in\FixedBoundary\Omega$ that see each other in~$\Omega$ inequality~\textup{\eqref{LightBound}} holds true. Suppose also that~$|\nabla L_x|+ |L_{x}(0)|$ is uniformly bounded when~$x$ belongs to a compact set.
If the function~$G$ is given by formula~\textup{\eqref{LightFunction},} then there exists some domain~$\Omega^*$ satisfying conditions~\textup{\eqref{FirstCondition}, \eqref{ThirdCondition}} such that~$\Omega$ is an extension of~$\Omega^*$ and~$G$ is locally concave on~$\Omega^*$.
\end{St}
\begin{proof}
Using Lemma~\ref{ConcavityNearTheBoundary}, we match each compact set~$K \subset \FixedBoundary \Omega$ its parameter~$\eps_K$ and define the function~$\rho\colon\FixedBoundary \Omega \to (0,1]$ by the formula
\begin{equation*}
\rho(x) = \min\Big(1,\sup\Big\{\frac {\eps_K}{2}\,\Big|\,\, K \,\,\hbox{is a compact subset of}\,\, \FixedBoundary\Omega,\,\, x \in K\Big\}\Big).
\end{equation*}
We build~$\Omega^*_1$ with the help of Lemma~\ref{FirstExtension} and this function~$\rho$. First, the domain~$\Omega$ is an extension of~$\Omega^*$. Second, if~$z \in \Omega^*$, then there exist two points~$x$ and~$y$ belonging to~$\FixedBoundary\Omega^*$ such that~$z \in [x,y] \subset \Omega^*$. There exists a compact set~$K \subset \FixedBoundary\Omega^*$ such that~$x \in K$ and~$\rho(x) \leq \eps_K$. In such a case,~$|x-y| \leq \eps_K$. By Lemma~\ref{ConcavityNearTheBoundary}, there exists a linear function~$L[G;z]$ such that~$G(z) = L[G;z](z)$ and~$G \leq L[G;z]$ in some neighborhood of~$z$. Therefore, the function~$G$ satisfies the assumptions of Fact~\ref{SuperDifferentialConcavity} in the domain~$\Omega^*$, thus it is locally concave.
\end{proof}

\paragraph{Proof of Theorem~\ref{ExtensionTheorem}.} We start with the function~$\BG_{\Omega,f}$ on the domain~$\Omega$. Fix~$\delta > 0$. Using Lemma~\ref{SmallStrictlyConcaveFunction} and Corollary~\ref{StrictlyConcaveUpperApproximation}, we build the strictly locally concave function~$\BG_{\Omega,f} + E_{\delta}$ such that~$\BG_{\Omega,f} \leq \BG_{\Omega,f} + E_{\delta} \leq \BG_{\Omega,f} + \delta$, and the family of linear functions~$\{L[\BG_{\Omega,f} + E_{\delta};x]\}_{x \in \FreeBoundary \Omega}$. Treating the curve~$\FreeBoundary \Omega$ as the fixed boundary of a searched-for domain and using Corollary~\ref{StrictlyConcaveMajorant} (with~$\Omega_1$ as~$\Omega_0$ in the conditions of the corollary and the functions~$L[\BG_{\Omega,f} + E_{\delta};x]$ as~$L_x$), we build some domain~$\hat{\Omega}'_{\delta}$ such that~$\FreeBoundary \Omega = \FixedBoundary \hat{\Omega}'_{\delta}$ and the ``boundary conditions''~$L[\BG_{\Omega,f} + E_{\delta};x]$ satisfy the concavity condition~\eqref{LightBound} ``inside'' this domain. In such a case, we construct the function~$\hat{G}_{\delta}$ by formula~\eqref{LightFunction}, i.e. by the formula
\begin{equation*}
\hat{G}_{\delta}(z) = \inf\{L[\BG_{\Omega,f} + E_{\delta};x](z)\big|\,\,[z,x] \subset \hat{\Omega}'_{\delta},\,\,x \in \FreeBoundary\Omega\},\quad z \in \hat{\Omega}'_{\delta}
\end{equation*} 
(the domain~$\hat{\Omega}'_{\delta}$ is needed for this construction). By Proposition~\ref{LightExtension}, the function~$\hat{G}_{\delta}$ is locally concave on some domain~$\hat{\Omega}_{\delta}$,~$\FixedBoundary\hat{\Omega}_{\delta} = \FreeBoundary\Omega$, satisfies the inequality~$\hat{G}_{\delta}(z) \leq L[\BG_{\Omega,f} + E_{\delta};x](z)$ in a neighborhood of~$x$ (for any~$x \in \FreeBoundary \Omega$),~$z \in \hat{\Omega}_{\delta}$, and the equality~$\hat{G}_{\delta}(x) = L[\BG_{\Omega,f} + E_{\delta};x](x) = \BG_{\Omega,f}(x) + E_{\delta}(x)$ for any~$x \in \FreeBoundary \Omega$. We define~$\tilde{G}_{\delta}$ to be equal~$\BG_{\Omega,f} + E_{\delta}$ on~$\Omega$ and~$\hat{G}_{\delta}$ on~$\hat{\Omega}_{\delta}$. The third point of Corollary~\ref{StrictlyConcaveUpperApproximation} implies that~$\tilde{G}_{\delta}(z) \leq L[\BG_{\Omega,f} + E_{\delta};x](z)$ in a neighborhood of~$x$ for any~$x \in \FreeBoundary\Omega$,~$z\in \tilde{\Omega}_{\delta}$. The function~$\tilde{G}_{\delta}$ is locally concave on~$\tilde{\Omega}_{\delta} = \Omega\cup\hat{\Omega}_{\delta}$ by Fact~\ref{SuperDifferentialConcavity} and satisfies the inequality~$\tilde{G}_{\delta} \leq \BG_{\Omega,f} + \delta$ on~$\Omega$ because the function~$\BG_{\Omega,f} + E_{\delta}$ does. \qed
\subsection{Proof of the main theorem for smooth case}\label{s45}
\begin{Th}\label{SmoothMainTheorem}
Suppose that the domain~$\Omega$ satisfies the assumptions~\textup{\eqref{FirstCondition},~\eqref{SecondCondition},~\eqref{ThirdCondition}} and also has~$C^2$-smooth fixed boundary. Let the function~$f$ be bounded from below and continuous. Then~$\BG_{\Omega,f} = \Bell_{\Omega,f}$.
\end{Th}
\begin{proof}
By Proposition~\ref{ContinuityOnFixedFoundary},~$\BG_{\Omega,f}$ is continuous at every point of the fixed boundary. Consequently, Theorem~\ref{GeometricMartingale} asserts that~$\BG_{\Omega,f} = \BM_{\Omega,f}$. By Corollary~\ref{MartingaleLessAnalytic},~$\BG_{\Omega,f} = \BM_{\Omega,f} \leq \Bell_{\Omega,f}$. By Corollary~\ref{AnalyticLessMartingale},~$\Bell_{\Omega,f} \leq \inf_{\tilde{\Omega}} \BM_{\tilde{\Omega},f} = \inf_{\tilde{\Omega}} \BG_{\tilde{\Omega},f}$; here the infimum is taken over all the extensions of~$\Omega$. By Corollary~\ref{ExtensionForSmoothFunctions},~$\inf_{\tilde{\Omega}} \BG_{\tilde{\Omega},f} = \BG_{\Omega,f}$, which finishes the proof.
\end{proof}
\begin{Th}\label{SmoothMainTheoremUnboundedFunction}
Suppose that the domain~$\Omega$ satisfies the assumptions~\textup{\eqref{FirstCondition},~\eqref{SecondCondition},~\eqref{ThirdCondition}} and also has~$C^2$-smooth fixed boundary. Let the function~$f$ be continuous. Then~$\BG_{\Omega,f} = \Bellb_{\Omega,f}$.
\end{Th}
The proof is the same.

The most technically difficult counterpart, Theorem~\ref{ExtensionTheorem}, is not needed for the classical domains~$\Omega$ (such as the parabolic strip or the hyperbolic strip; namely, \emph{all} the domains listed in Subsection~\ref{s13} on which the Bellman functions are considered) if one also assumes some mild summability conditions, because these domains have additional homogeneous structure. For example, see~\cite{IOSVZ}, proof of Statement~$2.6$, where the general extension theorem is replaced by a one line formula using homogeneity (however, the reasoning there requires additional summability, for example, with this method one is, seemingly, not able to prove Theorem~\ref{SuperSummability} below). It would be great either to find a shorter and more transparent proof of Theorem~\ref{ExtensionTheorem} in the full generality, or change it for something weaker that is still sufficient for the proof of the main theorem.

\section{Smoothing}\label{s5}
We begin this section with the construction that will be used as the main engine of the smoothing procedure. The idea resembles the one used to prove Theorem~\ref{ExtensionTheorem}, namely, we are going to define the smoothed function~$G'$ as an infimum of linear functions generated by the initial function~$G$, as we did it in formula~\eqref{LightFunction}. 

\begin{figure}[h!]
\includegraphics[height=9.5cm]{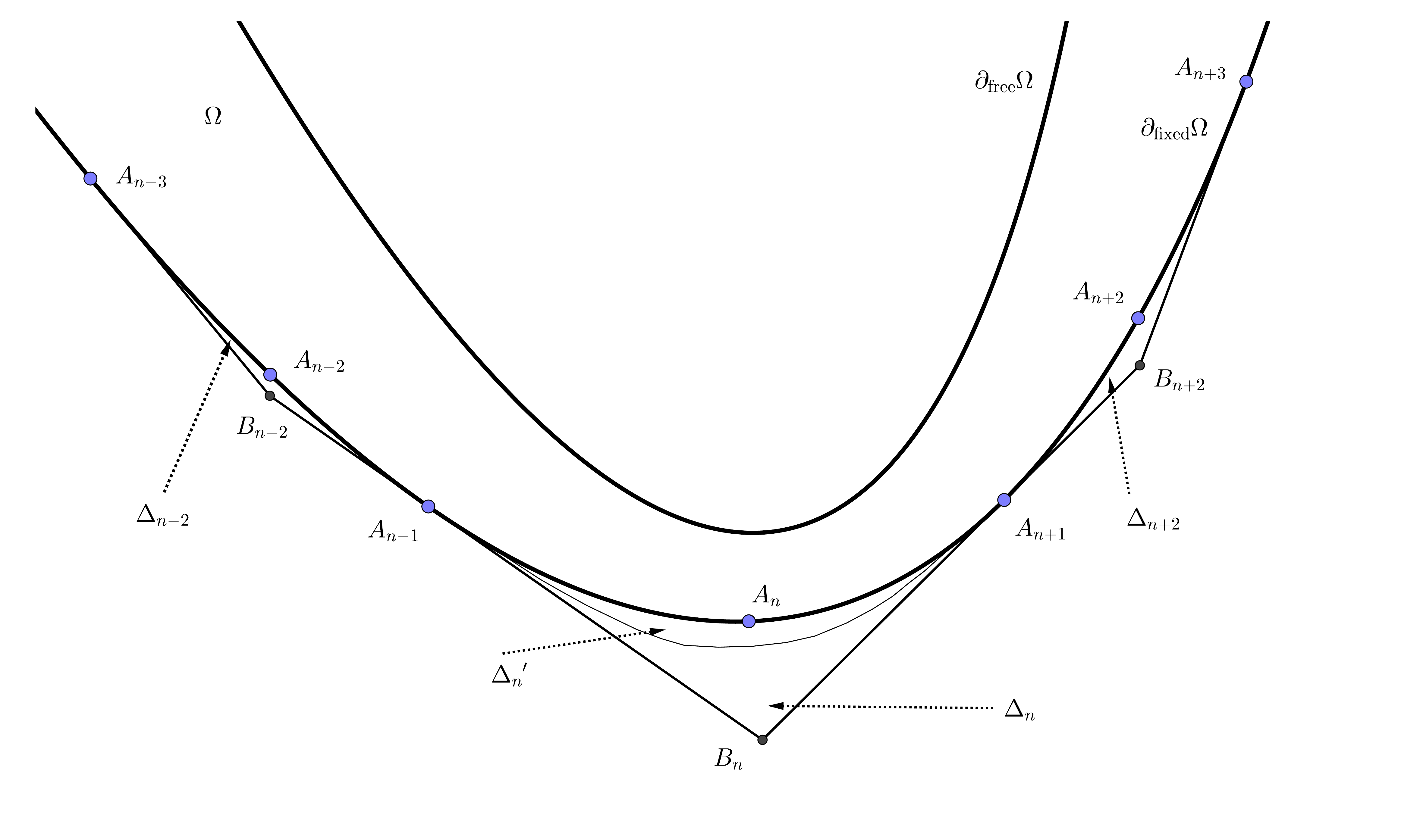}
\caption{Construction related to sequence~$\{A_n\}$.}
\label{fig:seq}
\end{figure}


Let~$\Omega$ satisfy conditions~\eqref{FirstCondition} and~\eqref{ThirdCondition}. Let~$G$ be a locally concave function on~$\Omega$ that is locally bounded. We choose a linear function~$L[G;x]$ for any~$x \in \interior \Omega$; as usual, we require two properties:~$L[G;x](x) = G(x)$ and~$G(y) \leq L[G;x](y)$ for any~$y$ such that~$[y,x] \subset \Omega$.

Let~$\{A_n\}_{n \in \mathbb{Z}}$ be a monotone sequence (the ordering was defined at the beginning of Subsection~\ref{s32}) of points on the fixed boundary such that~$[A_n,A_{n+2}]$ does not intersect~$\Omega_{1}$ and each point of~$\FixedBoundary \Omega$ lies between two points of this sequence (i.e. the sequence ``covers'' the whole boundary). We also require that at each point~$A_n$ the curve~$\FixedBoundary \Omega$ is differentiable. For each~$n$ we find a unique point~$B_n \in \mathbb{R}^2$ such that the segments~$[A_{n-1},B_n]$ and~$[B_n,A_{n+1}]$ touch~$\FixedBoundary\Omega$. The symbol~$\Delta_n$ denotes the open triangle with the vertices~$A_{n-1}$,~$B_n$, and~$A_{n+1}$.  We will also need one more domain that lies inside~$\Delta_n$. Let us call it~$\Delta_n'$. We introduce it in several steps. By~$H^{\lambda}_x$ we denote the homothety with coefficient~$\lambda$ and center~$x$ (i.e. the mapping~$\mathbb{R}^2 \ni y \mapsto x+\lambda(y-x)$). Consider the convex set
\begin{equation*}
H^{2}_{A_{n-1}}\big( \Omega \cap \Delta_n\big) \cap H^{2}_{A_{n+1}}\big( \Omega \cap \Delta_n\big).
\end{equation*}
The convex set~$\Delta_n'$ lies inside this set, contains~$\Omega \cap \Delta_n$, and has a~$C^2$-smooth boundary (except, maybe, at the points~$A_{n-1}$ and~$A_{n+1}$). See Figure~\ref{fig:seq} for the visualization of the construction. There are many methods how to construct such a set, we leave it to the reader. The main property that we will use is the following one: there is a small~$\eps$ such that when~$x \in \Delta_n' \setminus \Omega$ and~$|x - A_{n-1}| < \eps$, one has~$|x-x'| \leq |x' - A_{n-1}|$, where~$x' = [A_{n-1},x]\cap \FixedBoundary\Omega$; a similar condition also holds when~$A_{n+1}$ is replaced by~$A_{n-1}$.  

We are going to extend the function~$G$ through the fixed boundary. The following lemma is about extending through a small part of it.
\begin{Le}\label{LocalExtensionFixedBoundary}
In the above notation\textup, the function~$U_n\colon \mathbb{R}^2 \to \mathbb{R} \cup\{-\infty\}$ defined by the formula
\begin{equation}\label{UnFormula}
U_n(y) = \inf\{L[G;x](y)\mid x \in (A_{n-1}, A_{n+1})\}
\end{equation}
is finite inside~$\Delta_n$. If the function~$G$ is continuous at~$A_{n-1}$\textup, then the restriction of~$U_n$ to the domain~$\Delta'_n$ is also continuous at~$A_{n-1}$\textup(and the same claim holds with~$A_{n+1}$ in place of~$A_{n-1}$\textup).
\end{Le}
\begin{proof}
To prove the first claim, we fix a point~$y \in \Delta_n \setminus \Omega$ and introduce several points related to it. Let~$y_{n-1} = [A_{n-1},y]\cap \FixedBoundary\Omega$ and let~$y_{n+1} = [A_{n+1},y]\cap \FixedBoundary\Omega$. If~$x \in (A_{n-1},A_{n+1})$, then (the point~$x_y$ stands for~$[x,y]\cap (y_{n-1},y_{n+1})$, see Figure~\ref{fig:est})
\begin{equation}\label{LowerBound}
L[G;x](y) = L[G;x](x) + \frac{|y-x|}{|x_y - x|}\big(L[G;x](x_y) - L[G;x](x)\big) \geq -\frac{|y-x_y|}{|x_y - x|}G(x) + \frac{|y-x|}{|x_y - x|}G(x_y).
\end{equation}
This value is uniformly bounded from below when~$x \in (A_{n-1},A_{n+1})$, because the values~$G(x)$,~$G(x_y)$, and the coefficients are bounded. Therefore, the function~$U_n$ is finite on~$\Delta_n$, and Fact~\ref{ContinuityInInnerPoints} shows that it is continuous there.

\begin{figure}[h!]
\includegraphics[height=9.5cm]{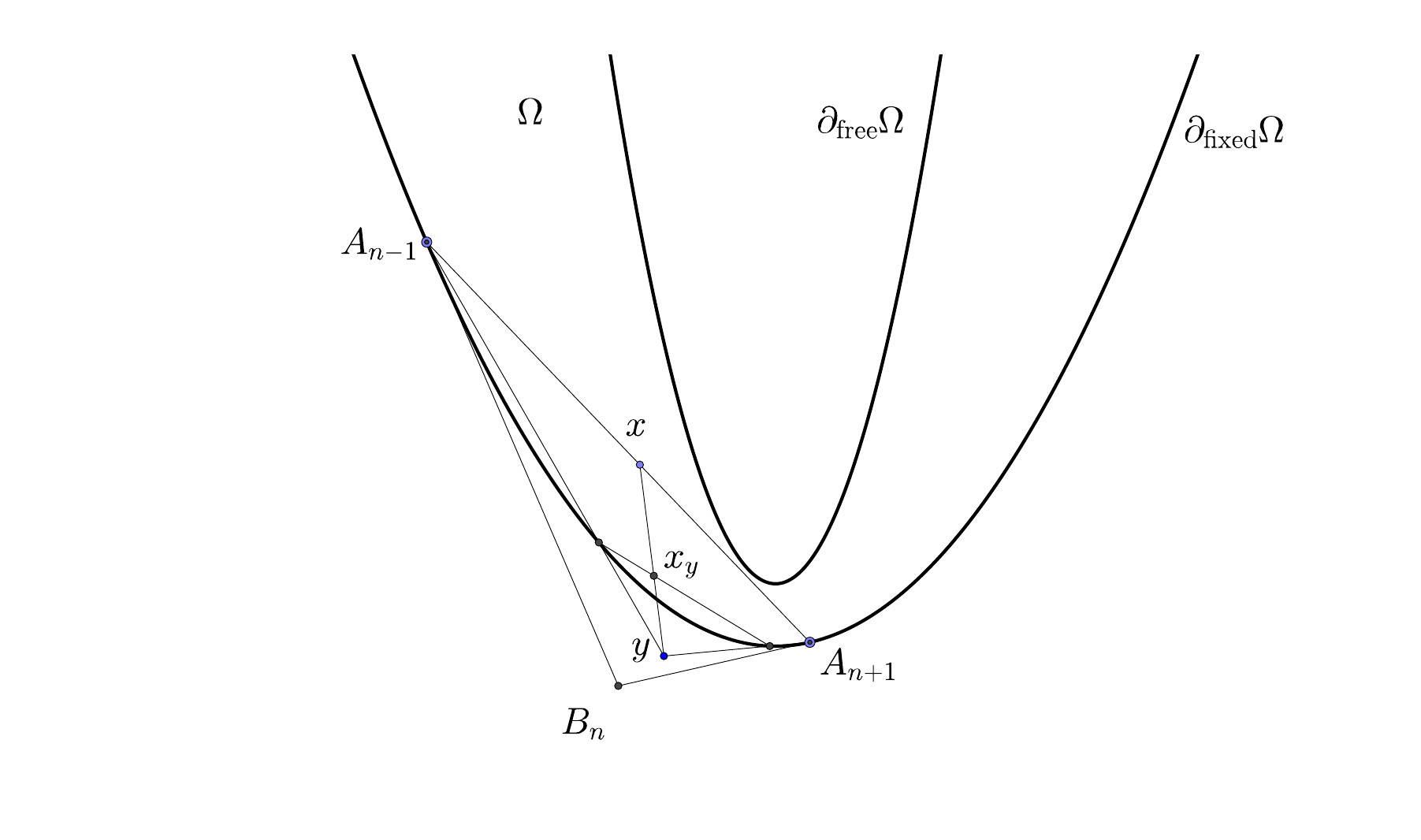}
\caption{Construction for point~$y$.}
\label{fig:est}
\end{figure}

Let us turn to the second claim of the lemma. First we prove the semi-continuity from below. Fix~$\eps > 0$. It suffices to prove a uniform in~$x \in (A_{n-1},A_{n+1})$ bound~$L[G;x](y) \geq G(A_{n-1}) - \eps$ when~$y \in \Delta_n'$ and~$|y - A_{n-1}|$ is sufficiently small. This estimate is obvious for~$y \in \Omega \cap \Delta_n$ due to the continuity of~$G$ at~$A_{n-1}$. In what follows, we suppose that~$y \in \Delta'_n \setminus \Omega$. Let~$r$ be a positive number such that~$\big|G(z) - G(A_{n-1})\big| < \delta$ when~$|z - A_{n-1}| < r$ and~$z \in \Omega$. Consider the two cases:~$|x - A_{n-1}| < r$ and~$|x-A_{n-1}| \geq r$.

In the case~$|x-A_{n-1}| \geq r$, we argue as follows: take~$|y-A_{n-1}|$ to be so small that~$\frac{|y-x_y|}{|x_y - x|} < \frac{1}{3}\eps(\sup_{[A_{n-1},A_{n+1}]} |G|)^{-1}$
whenever~$|x - A_{n-1}| \geq r$ and~$x \in (A_{n-1},A_{n+1})$ (one can do so because the curve~$\FixedBoundary\Omega$ is differentiable at the point~$A_{n-1}$); we also require that in such a case~$|x_{y} -A_{n-1}| < r$. Thus, formula~\eqref{LowerBound} provides the estimate
\begin{equation*}
L[G;x](y) \geq -\frac{\eps}{3} + \Big(1- \frac{\eps}{3}(\!\!\sup_{[A_{n-1},A_{n+1}]}\!\!\!\!\!\!\!\!|G|)^{-1}\Big)(G(A_{n-1}) - \delta) \geq G(A_{n-1}) - \eps \quad\hbox{if~$\delta$ is sufficiently small}.
\end{equation*}

In the case~$|x - A_{n-1}| \leq r$, we use the fact that the ratio~$\frac{|y-x_y|}{|x_y - x|}$ is bounded by~$2$. Indeed, if~$y$ is sufficiently close to~$A_{n-1}$, then~$\frac{|y-x_y|}{|x_y - x|} \leq \frac{|y - y_{n-1}|}{|y_{n-1} - A_{n-1}|}$ (it follows from the Thales theorem), and the latter quantity does not exceed~$2$ by our assumptions on the domain~$\Delta'_n$. Therefore, formula~\eqref{LowerBound} results in
\begin{equation*}
\begin{aligned}
L[G;x](y) \geq -\frac{|y-x_y|}{|x_y - x|}\big(G(A_{n-1}) +\delta\big) + \frac{|y-x|}{|x_y - x|}\big(G(A_{n-1}) -\delta\big) =\\ G(A_{n-1}) - \Big(2\frac{|y-x_y|}{|x_y - x|}+1\Big)\delta \geq G(A_{n-1}) -\eps 
\end{aligned}
\end{equation*}
if~$\delta \leq \frac{\eps}{7}$.
The semi-continuity from below  of~$U_n$ at~$A_{n-1}$ is proved.

 On the other hand, ~$G$ is continuous at the point~$A_{n-1}$ (and concave), therefore~$G(A_{n-1})$ equals~$\lim_{x\to A_{n-1}}L[G;x](A_{n-1})$, where~$x \in [A_{n-1},A_{n+1}]$. Fix~$\eps > 0$ and let~$x\in [A_{n-1},A_{n+1}]$ be such that~$G(A_{n-1}) \geq L[G;x](A_{n-1}) -\eps$. Then,
\begin{equation*}
\varlimsup\limits_{y \to A_{n-1}} U_n(y) \leq \varlimsup\limits_{y \to A_{n-1}} L[G;x](y) = L[G;x](A_{n-1}) \leq G(A_{n-1}) + \eps.
\end{equation*}
Sending~$\eps$ to~$0$, we get the semi-continuity from above.
\end{proof}

\begin{Le}\label{ExtensionThroughFixedBoundaryLemma}
 Let~$G$ be a locally concave function on the domain~$\Omega$ satisfying conditions~\textup{\eqref{FirstCondition}, \eqref{SecondCondition}}. There exists a locally concave function~$G'$ on the domain~$\Omega \cup (\cup_{k\in\mathbb{Z}}\Delta'_{2k})$ such that~$G = G'$ on~$\Omega \setminus(\cup_{k\in\mathbb{Z}}\Delta'_{2k})$ and~$G\leq G'$ on~$\Omega$\textup, and~$G'$ is continuous everywhere\textup, except\textup, possibly\textup, the points~$A_{2k-1}, k \in \mathbb{Z}$. If~$G$ is continuous at~$A_{2k-1}$ for some numbers~$k$\textup, then~$G'$ can be chosen to be continuous at the same points.  
\end{Le} 
\begin{proof}
For each~$n \in \mathbb{Z}$, define the function~$U_n\colon \mathbb{R}^2 \to \mathbb{R}$ by formula~\eqref{UnFormula}.
This function is concave, it is not less than~$G$ on~$\Omega \cap \cl \Delta_n$, and coincides with~$G$ on the segment~$(A_{n-1},A_{n+1})$.
 
Let the function~$G'$ be equal~$G$ on~$\Omega \setminus \cup_k\Delta'_{2k}$ and~$U_{2k}$ on~$\Delta'_{2k}$ for each~$k \in \mathbb{Z}$. First, this function is equal to~$G$ on~$\Omega \setminus\cup_{k}\Delta'_{2k}$. Second, this function is locally concave by Fact~\ref{SuperDifferentialConcavity}. Third, it is continuous at every point of its definition domain (because the functions~$U_{2k}$ are finite on the corresponding domains by Lemma~\ref{LocalExtensionFixedBoundary}), except for the points~$A_{2k-1}$,~$k\in\mathbb{Z}$, by Fact~\ref{ContinuityInInnerPoints}. The continuity result for the points~$A_{2k-1}$ has also been proved in Lemma~\ref{LocalExtensionFixedBoundary}.
\end{proof}
\begin{Cor}\label{NonSmoothBoundaryCor}
Let~$\Omega$ satisfy conditions~\textup{\eqref{FirstCondition}, \eqref{SecondCondition}, \eqref{ThirdCondition}}. If the function~$f\colon\FixedBoundary \to \mathbb{R}$ is bounded from below and continuous\textup, then~$\Bell_{\Omega,f} = \BG_{\Omega,f}$. If~$f$ is any continuous function\textup, then~$\Bellb_{\Omega,f} = \BG_{\Omega,f}$.
\end{Cor}
\begin{proof}
We prove only the first assertion, the second one is similar. There is nothing to prove if~$\BG$ is infinite, so we assume that it is locally bounded. The inequality~$\BG_{\Omega,f} \leq \Bell_{\Omega,f}$ follows from Theorem~\ref{GeometricMartingale} and Corollary~\ref{MartingaleLessAnalytic}. 

To prove the converse inequality, we fix some point~$x \in \Omega\setminus\FixedBoundary\Omega$.
Fix~$\eps > 0$ and choose a monotone sequence~$A_k$ of the points in~$\FixedBoundary \Omega$ such that~$x \notin \Delta_k$ for any~$k$ and the values of~$f$ inside the arc between~$A_{k-3}$ and~$A_{k+3}$ do not differ for more than~$\eps$ (so, on this arc~$f$ is a constant~$f_k$ plus an error term that does not exceed~$\eps$ in modulus). We apply Lemma~\ref{ExtensionThroughFixedBoundaryLemma} to the locally concave function~$\BG_{\Omega,f}$ to get the function~$G'$. This function is locally concave and continuous on the domain~$\Omega \cup (\cup_{k\in\mathbb{Z}}\Delta'_{2k})$, moreover,~$G' \geq f_{2k} - \eps$ on the arc of~$\FixedBoundary\Omega$ between the points~$A_{2k-1}$ and~$A_{2k+1}$. We can construct a convex closed set~$\Omega_0'$ such that~$\Omega_0 \subset \Omega_0' \subset \Omega_0 \cup (\cup_k \Delta'_{2k})$,~$G'\big|_{\partial\Omega_0'\cap\Delta_{2k}} \geq f_{2k} -2\eps$, and the boundary of~$\Omega_0'$ is~$C^2$-smooth except, maybe, for the points~$A_{2k-1}$,~$k \in \mathbb{Z}$.  On the fixed boundary of this new domain, we find a new sequence (say,~$\{A_k'\}_k$) such that its points differ from the~$A_k$ (and between each points~$A_k$ and~$A_{k+1}$ there is at least two points of the new sequence; we may also suppose that the tangent triangles that correspond to this sequence lie inside~$\Omega \cup \cup_n \Delta_n$), but the point~$x$ does not belong to the corresponding triangles~$\Delta_k''$. Applying Lemma~\ref{ExtensionThroughFixedBoundaryLemma} to the domain~$\Omega' = \cl\Omega_0' \setminus \Omega_1$ and the function~$G'$ (with the sequence~$\{A_k'\}_k$), and then narrowing the obtained domain near the points~$A_{2k}'$ (where its boundary may be non-smooth), we get the domain~$\Omega_0''$ with a~$C^2$-smooth boundary such that~$\Omega \subset \Omega_0''$, and a locally concave function~$G''$,~$G''$ is continuous and locally concave on~$\cl \Omega_0'' \setminus \Omega_1$, and~$G''(x) = \BG_{\Omega,f}(x)$. We may also require that~$G'' \geq f_k - 4\eps$ on~$\FixedBoundary \Omega'' \cap \Delta_k$ (by one more narrowing). Let~$\Omega'' = \cl\Omega_0'' \setminus \Omega_1$, then
\begin{equation*}
\BG_{\Omega'',G''|_{\FixedBoundary \Omega''}}(z) \geq f_k - 4\eps \geq f(z) - 5\eps,\quad  z\in\FixedBoundary\Omega \cap \Delta_{k}.
\end{equation*}
 We apply Corollary~\ref{ExtensionForSmoothFunctions} to the function~$\BG_{\Omega'',G''|_{\FixedBoundary\Omega''}}$ on the domain~$\Omega''$ to find an extension~$\tilde{\Omega}''$ of the domain~$\Omega''$ such that~$\BG_{\tilde{\Omega}'',G''|_{\FixedBoundary\Omega''}}$ does not exceed~$\BG_{\Omega'',G''|_{\FixedBoundary\Omega''}} + \eps$ on the domain~$\Omega''$. Let~$\tilde{\Omega} = \tilde{\Omega}'' \setminus(\Omega'' \setminus \Omega)$ and let~$\tilde{f}$ be the restriction of~$\BG_{\tilde{\Omega}'',G''|_{\FixedBoundary\Omega''}} + 5\eps$ to~$\FixedBoundary 
 \Omega$. Then,~$\tilde{f} \geq f$ and~$\BG_{\tilde{\Omega},\tilde{f}}$ does not exceed~$\BG_{\Omega,f} + 6\eps$ at the point~$x$. We can write
\begin{equation*}
\Bell_{\Omega,f}(x)\!\!\!\!\stackrel{\scriptscriptstyle \hbox{\tiny Cor.~}\ref{AnalyticLessMartingale}}{\leq}\!\!\!\! \BM_{\tilde{\Omega},f}(x) \stackrel{\scriptscriptstyle \hbox{\tiny Th.~}\ref{GeometricMartingale}}{=} \BG_{\tilde{\Omega},f}(x) \leq \BG_{\tilde{\Omega},\tilde{f}}(x) \leq \BG_{\Omega,f}(x) + 6\eps.
\end{equation*}
Sending~$\eps$ to~$0$, we prove the corollary.
\end{proof}
Now we are ready to prove the main theorem in the full generality.
\begin{Cor}\label{DiscontinuousFunctionNonSmoothBoundary}
Let~$\Omega$ satisfy conditions~\textup{\eqref{FirstCondition}, \eqref{SecondCondition}, \eqref{ThirdCondition}}. If the function~$f\colon\FixedBoundary \to \mathbb{R}$ is bounded from below\textup, then~$\Bell_{\Omega,f} = \BG_{\Omega,f}$. If~$f$ is any locally bounded function\textup, then~$\Bellb_{\Omega,f} = \BG_{\Omega,f}$.
\end{Cor}
\begin{proof}
We will prove only the first assertion, the second one is similar. There is nothing to prove if~$\BG$ is infinite, so we assume that it is locally bounded. It suffices to prove that
\begin{equation*}
\av{f(\varphi)}{I} \leq \BG_{\Omega,f}(\av{\varphi}{I})
\end{equation*}
for any non-constant~$\varphi \in \Class_{\Omega}$.
Choose a monotone sequence~$A_k$ of the points in~$\FixedBoundary \Omega$ such that~$\av{\varphi}{I} \notin \Delta_k$ for any~$k$. Let~$G'$ denote the function provided by Lemma~\ref{ExtensionThroughFixedBoundaryLemma} applied to the locally concave function~$\BG_{\Omega,f}$. The function~$G'$ is continuous on~$\Omega$ except for the points~$A_{2k-1}$,~$k \in \mathbb{Z}$. Let~$G''$ denote the function which is obtained by an application of the same lemma to the function~$G'$ and the sequence of poins~$A'_k \df A_{k-1}$,~$k \in \mathbb{Z}$ (so the corresponding triangles have the points~$A_{2k}$ among their vertices). The function~$G''$ has several nice properties: it is continuous on the whole domain~$\Omega$ (the main point is that Lemma~\ref{ExtensionThroughFixedBoundaryLemma} allows it to remain continuous at the points~$A_{2k-1}$), it is not less than~$\BG_{\Omega,f}$, and
\begin{equation*}
G''(\av{\varphi}{I}) = \BG_{\Omega,f}(\av{\varphi}{I}).
\end{equation*}
The proof of the inequality wanted looks like this:
\begin{equation*}
\av{f(\varphi)}{I} \leq \av{G''(\varphi)}{I} \!\!\!\stackrel{\scriptscriptstyle\hbox{\tiny Cor.~}\ref{NonSmoothBoundaryCor}}{\leq} \BG_{\Omega,G''\big|_{\FixedBoundary\Omega}}\!\!\!\!\!\!\!\!\!\!\!\!(\av{\varphi}{I}) \leq G''((\av{\varphi}{I})) = \BG_{\Omega,f}(\av{\varphi}{I}).
\end{equation*}
\end{proof}

\section{Summability questions}\label{s6}
We have statements of two types: either the boundary value~$f$ is bounded from below or we consider functions (and martingales) with bounded values. In this section we try to combine them, i.e. to extract the \emph{a posteriori} correctness of the integral in formula~\eqref{BellmanFunction} from the finiteness of~$\BG_{\Omega,f}$.

We begin with a discouraging claim: there exists a function~$f$ such that the function~$\Bellb_{\Omega,f} = \BG_{\Omega,f}$ is finite, but for some~$\varphi \in \Class_{\Omega}$ the integral~$\av{f(\varphi)}{I}$ is not well defined (i.e. both~$f_+(\varphi)$ and~$f_-(\varphi)$ are not summable,~$f_+ = \min(f,0)$,~$f_- = \max(f,0)$). Such an example can be constructed for the case of the parabolic strip (thus the class~$\Class_{\Omega}$ is isomorphic to~$\BMO_{\eps}$ for some~$\eps$) using formulas from the paper~\cite{IOSVZ}. We do not dwell on this, the reader is referred to the forthcoming paper~\cite{ISVZ}. The example shows that the summablity questions are not only technical ones. First, we prove that the summability of~$f(\varphi)$ for all~$\varphi \in \Class_{\Omega}$ is equivalent to finiteness of~$\BG_{\Omega,f_+}$ (under assumption that~$\BG_{\Omega,f}$ is finite), and second, we provide an easy-to-verify criterion for the function~$f$ that gives the finiteness of~$\BG_{\Omega,f}$ and~$\BG_{\Omega,f_+}$ (and thus the correctness of formula~\eqref{BellmanFunction}).

In what follows, we say that the Lebesgue integral~$\int g$ is well defined if either~$\int g_-$ or~$\int g_+$ is finite.

\begin{St}\label{Summability}
Suppose that~$\Omega$ satisfies conditions~\textup{\eqref{FirstCondition}, \eqref{SecondCondition}, \eqref{ThirdCondition}}\textup, the function~$f$ is locally bounded from below\textup, and~$\BG_{\Omega,f}$ is finite. Then\textup, the integral~$\av{f(\varphi)}{I}$ is well defined for all~$\varphi \in \Class_{\Omega}$ if and only if~$\BG_{\Omega,f_+}$ is finite. Moreover\textup, in such a case~$\Bell_{\Omega,f} = \BG_{\Omega,f}$.
\end{St}
\begin{proof}
Suppose that~$\BG_{\Omega,f_+}$ is finite. The function~$f_+$ is bounded from below, therefore, by Corollary~\ref{DiscontinuousFunctionNonSmoothBoundary},~$\Bell_{\Omega,f_+} = \BG_{\Omega,f_+} < \infty$, which means that~$\av{f_+(\varphi)}{I}$ is always finite. Thus,~$\av{f(\varphi)}{I}$ is well defined.

The reverse implication is a bit more difficult. First, we prove that if~$\av{f(\varphi)}{I}$ is well defined for all~$\varphi \in \Class_{\Omega}$ and~$\BG_{\Omega,f}$ is finite, then~$\Bell_{\Omega,f} = \BG_{\Omega,f}$ (and thus~$\Bell_{\Omega,f}$ is finite). By Corollary~\ref{DiscontinuousFunctionNonSmoothBoundary}, we know that~$\BG_{\Omega,f} = \Bellb_{\Omega,f} \leq \Bell_{\Omega,f}$, so it suffices to prove the inequality
\begin{equation*}
\av{f(\varphi)}{I} \leq \BG_{\Omega,f}(\av{\varphi}{I})
\end{equation*}
for any~$\varphi \in \Class_{\Omega}$. We may assume that~$\varphi$ is non-constant. By Corollary~\ref{MonotonicRearrangementCor}, we may assume that~$\varphi$ is defined on~$[0,1]$ and monotonic. Let~$\varphi_n$ be the restriction of~$\varphi$ to the interval~$I_n = [\frac{1}{n},\frac{n-1}{n}]$,~$n > 2$. For each~$n$,~$\varphi_n \in \Classb_{\Omega}$, thus,
\begin{equation*}
\av{f(\varphi_n)}{I_n} \leq \BG_{\Omega,f}(\av{\varphi_n}{I_n}).
\end{equation*}
The integral~$\av{f(\varphi)}{[0,1]}$ is well defined, thus~$\av{f(\varphi_n)}{I_n} \to \av{f(\varphi)}{[0,1]}$ as~$n \to \infty$. By continuity of the function~$\BG$, provided by Fact~\ref{ContinuityInInnerPoints} and Proposition~\ref{ContinuityOnFreeBoundary},~$\BG_{\Omega,f}(\av{\varphi_n}{I_n}) \to \BG_{\Omega,f}(\av{\varphi}{[0,1]})$, and the inequality is proved. So,~$\Bell_{\Omega,f} = \BG_{\Omega,f}$. 

We have to prove that~$\Bell_{\Omega,f_+}$ is finite. Assume the contrary, let it be infinite at some point~$x$. By the very definition, there exists a sequence~$\{\varphi_n\}_n$ of functions belonging to~$\Class_{\Omega}$ such that~$\av{\varphi_n}{I} = x$ and~$\av{f_+(\varphi_n)}{I} \geq 2^n$. We apply Corollary~\ref{Concatenation} to this sequence of functions and the sequence of numbers~$\alpha_n =2^{-n}$ to obtain a function~$\varphi \in \Class_{\Omega}$ such that~$\av{\varphi}{I} = x$, but~$\av{f_+(\varphi)}{I} = \infty$. The integral~$\av{f(\varphi)}{I}$ is well defined, therefore,~$\av{f(\varphi)}{I} = \infty$, which contradicts the finiteness of~$\BG_{\Omega,f}$.
\end{proof}

Though we have a criterion for the Bellman function to be finite given by the previous proposition, it seems to be rather difficult to verify. Now we give a condition on~$f$ that guarantees the finiteness of~$\Bell_{\Omega,f_+}$. It is explicit modulo calculations of certain Bellman functions, which are very easy, but need some technique that requires a vast exposition. So, we write our estimates in terms of certain Bellman functions, and then show, what do they give for the case of the parabolic strip (in particular, we show that our estimates are not far from being sharp).

We parameterize the fixed boundary with the inverse projection as we did in the beginning of Subsection~\ref{s32}. Next, we choose an increasing sequence~$\{A_n\}_{n}$ of points lying on the fixed boundary such that for each~$n$ the segment~$[A_{n-1},A_n]$ touches the upper boundary (the sequence may be finite, for example, for the domain~$\{x\in\mathbb{R}^2\mid x_1x_2 \geq 1, (x_1-1)(x_2-1) \leq 2, x_1,x_2 > 0\}$). Let~$\mathfrak{A}_n$ be the arc of the fixed boundary with the endpoints~$A_{n-1}$ and~$A_{n+1}$. Define the maximal function~$f^{\flat}\colon \FixedBoundary\Omega\to \mathbb{R}_+$ by the formula
\begin{equation}\label{MF}
f^{\flat} = \sum\limits_{k\in\mathbb{Z}}\sup\big\{ |f(z)|\,\big|\,\, z \in \mathfrak{A}_k\big\} \chi_{_{\mathfrak{A}_k}}.
\end{equation}  
This maximal function approximates~$f$, if~$f$ does not oscillate to much. It suffices to find at least one finite function in the class~$\Lambda_{\Omega,f}$. For example, it can be the function
\begin{equation*}
\sum\limits_{k\in\mathbb{Z}} f^{\flat}(A_{k})\BG_{\Omega,\chi_{_{\mathfrak{A}_k}}}
\end{equation*}
if the series converges at every point of the domain (a sum of locally concave functions is locally concave). The Bellman functions on the right-hand side can be easily calculated. As we have said, we would not do this in the full generality, but concentrate on the example of the ball in~$\BMO$. In this case, one can take~$A_k = (2k,4k^2)$, and write an easy estimate for the corresponding Bellman function:
\begin{equation*}
\BG_{\Omega,\chi_{_{\mathfrak{A}_k}}}(x) \lesssim e^{-|x_1-2k|}.
\end{equation*}
Similar Bellman functions had been calculated in~\cite{Osekowski2} and~\cite{Vasyunin4}. Thus, we have the following theorem that lies beyond the classical John--Nirenberg inequality.
\begin{Th}\label{SuperSummability}
If the function~$f$ satisfies the condition
\begin{equation*}
\sum\limits_{k\in\mathbb{Z}} e^{-|k|/\eps}\sup\limits_{[k-2,k+2]} |f| < \infty,
\end{equation*}
for some~$\eps > 0$\textup, then~$\av{f(\varphi)}{I}$ is finite for any~$\varphi \in \BMO_{\eps}(I)$.
\end{Th}
This theorem is sharper than the condition~$|f(t)| \leq e^{|t|/\eps_+}$ for some~$\eps_+ > \eps$ from~\cite{IOSVZ}. However, Theorem~\ref{SuperSummability} heavily relies on Theorem~\ref{ExtensionTheorem}.

\section{Generalizations}\label{s7}
As the reader may have already noticed, the abstract theory of Subections~\ref{s22} and~\ref{s23} can be transferred to Riemannian manifolds (or even more general objects). This may be useful in the context of working with domains of the type ``circle minus circle'', e.g.~$\Omega = \{x \in \mathbb{R}^2 \mid \eps \leq x_1^2 + x_2^2 \leq 1\}$. For such a domain one can consider similar function classes, Bellman functions, etc.. However, one does not have here analogs of the embedding results from Subsections~\ref{s31} and~\ref{s32}. Moreover, the natural domain for the martingales is not the domain~$\Omega$ itself, but its fundamental covering (because the martingales do not ``feel'' the fundamental group). However, the corresponding class of functions on the fixed boundary of this manifold does not match the naturally defined one. Nevertheless, the main theorem still holds. This is again due to topological reasons. Namely, the martingales that optimize the Bellman function cannot turn over the inner circle. This is a consequence of the following conjecture.
\begin{Conj}
Let~$\Omega_0$ be a compact convex subset of~$\mathbb{R}^2$ with non-empty interior. Suppose that~$\Omega_1$ is a convex open set such that~$\cl\Omega_1 \subset \Omega_0$. Let~$\Omega =\Omega_0 \setminus \Omega_1$. In such a case\textup, for any bounded function~$f\colon \FixedBoundary \Omega \to \mathbb{R}$ there are at least two points~$x_1,x_2$ on~$\FreeBoundary\Omega$ such that~$\BG_{\Omega,f}(x_i) = \BG_{\Omega_0,f}(x_i)$\textup,~$i=1,2$.
\end{Conj} 
In particular, this shows that the foliation (see~\cite{IOSVZ} for the definition) consisting of the tangents only cannot be a foliation of some minimal locally concave function (except for the trivial case of a linear function). So, the study of minimal locally concave functions and corresponding problems for functions and martingales may be influenced by some topological properties of the domain (here the effect comes from the fact the domain is not simply connected). 

Conjecture~\ref{Manifold} (see Appendix) seems to be very natural. It is also required for further study of abstract locally concave functions. 

A more demanding thing is to find some similar duality for dyadic problems (see, e.g.~\cite{SV} for some of them). There  is surely some: one can consider dyadic martingales on domains. However, the corresponding properties of minimal dyadically concave functions are not studied in a similar way that is presented here. There are also some small miracles: for example, the solution of the John--Nirenberg problem for dyadic classes appears to coincide with the solution of continuous problem for some bigger ball (or what is the same, for some wider strip), see~\cite{SV,SV2}.

There is a conjecture that, in a sense, related to the Bellman function on dyadic classes.
\begin{Conj}
Define the Bellman function for the John--Nirenberg inequality on a cube by the formula
\begin{equation*}
\Bell_{\BMO,d,\eps,\exp}(x_1,x_2) = \sup\Big\{\av{\exp(\varphi)}{Q}\mid \av{\varphi}{Q} = x_1, \av{\varphi^2}{Q} = x_2, \|\varphi\|_{\BMO(Q)} \leq \eps\Big\}.
\end{equation*}
Here~$Q$ is a~$d$-dimensional cube \textup(and the~$\BMO$ norm is defined as the supremum of the value~\eqref{AlmostBMO}\textup, where~$J$ runs through all subcubes of~$Q$\textup).
There exists some~$c = c(d) > 1$ such that~$\Bell_{\BMO,d,\eps,\exp} =\Bell_{\BMO,1,c(d)\eps,\exp}|_{\Omega_{\eps}}$\textup, where~$\Omega_{\eps} = \{x \in \mathbb{R}^2\mid x_1^2\leq x_2 \leq x_1^2 + \eps^2\}$.
\end{Conj}
The conjecture says that the Bellman function for the~$d$-dimensional problem can be obtained by restricting  the solution of the corresponding one-dimensional problem on some bigger ball to~$\Omega_{\eps}$. The exponential function is important here, because it gives additional homogeneous structure.

As we have seen, this problem is a friend of the corresponding monotonic rearrangement estimate.
\begin{Conj}
The constant~$c(d)$ from the previous conjecture is the best possible constant in the inequality
\begin{equation*}
\|f^*\|_{\BMO([0,1])} \leq c(d)\|f\|_{\BMO(Q)}.
\end{equation*}
\end{Conj}

We also believe that considering zigzag martingales (see~\cite{Osekowski} for the definition) on domains and manifolds may lead to new estimates for the martingale transform, i.e. martingales on domain, or~$w$-martingales may build the bridge between the two very nearby, but still not very connected on the formal level, topics: the Bellman function and the Burkholder method.

\appendix\label{appendix}\appendixpage
\addappheadtotoc
\section{Projective transform trick}\label{SPTT}
We describe a procedure that enables us to work with the domains~$\Omega$ for which~$\cl\Omega_0\notin \Omega_1$ (for example, with the natural domain for the Gehring class). We are not trying to reach  the full generality, but give a recipe that allows to cover all the known cases. Let~$\Phi$ and~$\Psi$ be strictly convex non-negative functions on the half-line~$[0,\infty)$. Let them grow faster than any linear function at infinity and let~$\Phi'(0) = \Psi'(0) = 0$. Suppose also that~$\Phi < \Psi$ on the open half-line. Consider the sets~$\Omega_0 = \{x \in \mathbb{R}^2 \mid x_1 > 0, x_2 > \Phi(x_1)\}$ and~$\Omega_1 = \{x \in \mathbb{R}^2 \mid x_1 > 0, x_2 > \Psi(x_1)\}$. As usual, 
\begin{equation*}
\Omega = \big\{x \in \mathbb{R}^2\,\big|\,\, x_1 > 0, \,\,\Phi(x_1) \leq x_2 \leq \Psi(x_1)\big\}.
\end{equation*}
The subset of~$\Omega$ where the inequality~$\Phi(x_1) = x_2$ turns into equality is called the fixed boundary, the subset of~$\Omega$ where the inequality~$x_2 = \Psi(x_1)$ turns into equality, is called the free boundary. As usual, we can define the class~$\Class_{\Omega}$ by formula~\eqref{AnalyticClass} and the Bellman function by formula~\eqref{BellmanFunction}. The dual problem is also well defined by formulas~\eqref{Lambdaclass} and~\eqref{MinimalLocallyConcave}. 

We are going to perform a projective transform. Let~$\omega$ be a subset of~$\mathbb{R}^3$ given by the formula
\begin{equation*}
\omega = \{z \in \mathbb{R}^3\mid z = (tx_1,tx_2,t), \,\,t > 0, x \in \Omega\}.
\end{equation*}
For any function~$G \colon \Omega \to \mathbb{R} \cup \{+\infty\}$ define the function~$g\colon \omega \to \mathbb{R} \cup \{+\infty\}$ by the formula~$g(tx,t) = tG(x)$. It is not hard to see that the function~$g$ is locally concave on~$\omega$ if and only if~$G$ is locally concave on~$\Omega$ (for details, see~\cite{ISZ}, Subsection~$2.3$). Consider now the section of the three-dimensional problem by the hyperplane~$\{z \in \mathbb{R}^3\mid z_1 = 1\}$. Let the coordinates there be~$y_1 = z_3$ and~$y_2 = z_2$. Then, the section of the domain~$\omega$ looks like this:
\begin{equation*}
\Omega_{\pr} =\big\{y \in \mathbb{R}^2\,\big|\,\, y_1 > 0, \,\,\Phi_{\pr}(y_1) \leq y_2 \leq \Psi_{\pr}(y_1)\big\}, \quad \hbox{where}\,\,\Phi_{\pr}(s) = s\Phi(s^{-1}), \Psi_{\pr}(s) = s\Psi(s^{-1}).
\end{equation*}
On the level of points, our projective transform acts by the rule~$(x_1,x_2)\mapsto (\frac{1}{x_1},\frac{x_2}{x_1}) = (y_1,y_2)$. It is not hard to see that the domain~$\Omega_{\pr}$ satisfies the assumptions~\eqref{FirstCondition} and~\eqref{ThirdCondition}. The assumption~\eqref{SecondCondition} is equivalent to the~$C^2$-smoothness of the function~$\Psi$. For any~$G \in \Lambda_{\Omega,f}$ we consider the function~$G_{\pr}$ given by the formula~$G_{\pr}(y_1,y_2) = \frac{G(x_1,x_2)}{x_1} = y_1G(\frac{1}{y_1},\frac{y_2}{y_1})$. Again,~$G_{\pr}$ is the restriction of the homogeneous function~$g$ to the hyperplane~$\{z \in \mathbb{R}^3\mid z_1 = 1\}$, therefore,~$G_{\pr}$ is locally concave whenever~$g$ is locally concave. Thus,~$G \in \Lambda_{\Omega,f}$ whenever~$G_{\pr} \in \Lambda_{\Omega_{\pr},f_{\pr}}$, where~$f_{\pr}\colon \FixedBoundary \Omega_{\pr}\to \mathbb{R}$ is given by the formula~$f_{\pr}(y_1,y_2) = \frac{f(x_1,x_2)}{x_1}$. So, we see that
\begin{equation*}
\BG_{\Omega_{\pr},f_{\pr}}(y_1,y_2) = \frac{\BG_{\Omega_{\pr},f_{\pr}}(x_1,x_2)}{x_1},\quad y \in \Omega_{\pr}.
\end{equation*}

The project transform acts on~$\Class_{\Omega}$ in a more sophisticated way. For any function~$\varphi \in \Class_{\Omega}$ denote its first coordinate by~$\varphi_1$ (thus,~$\varphi(t) = (\varphi_1,\Phi(\varphi_1))$). For convenience, we assume that~$I = [0,1]$. We note that the function~$t \mapsto \av{\varphi_1}{[0,1]}^{-1}\int_0^t \varphi_1(\zeta)\,d\zeta$ is an increasing absolutely continuous function mapping~$[0,1]$ to~$[0,1]$. Let~$\psi$ be its inverse function. Denote the projective transform of~$\varphi_1$ by the formula
\begin{equation*}
(\varphi_1)_{\pr}(t) = \frac{1}{\varphi_1(\psi(t))}.
\end{equation*}
We claim that the corresponding function~$\varphi_{\pr} = ((\varphi_1)_{\pr},\Phi_{\pr}((\varphi_1)_{\pr}))$ belongs to the class~$\Class_{\Omega_{\pr}}$ (and vice versa, if~$\varphi_{\pr} \in \Class_{\Omega_{\pr}}$, then~$\varphi \in \Class_{\Omega}$). Indeed, it suffices to prove that for any interval~$J_{\pr} \subset [0,1]$ there exists an interval~$J \subset [0,1]$ such that~$\av{\varphi_{\pr}}{J_{\pr}} = (\av{\varphi}{J})_{\pr}$. But this is straightforward (for clarity, let~$J_{\pr} = [a,b]$):
\begin{equation*}
\av{\varphi_{\pr}}{J_{\pr}} = \Bigg(\frac{\int_a^b \frac{ds}{\varphi_1(\psi(s))}}{\int_a^b ds},\frac{\int_a^b \frac{\Phi\big(\varphi_1(\psi(s))\big)ds}{\varphi_1(\psi(s))}}{\int_a^b ds}\Bigg) \stackrel{\scriptstyle \psi(s) =t}{=}\Big(\frac{1}{\av{\varphi_1}{[\psi(a),\psi(b)]}}, \frac{\av{\Phi(\varphi_1)}{[\psi(a),\psi(b)]}}{\av{\varphi_1}{[\psi(a),\psi(b)]}}\Big) = (\av{\varphi_1}{J},\av{\Phi(\varphi_1)}{J})_{\pr}, 
\end{equation*}
where~$J = [\psi(a),\psi(b)]$. The same change of a variable justifies the equality~$\av{f(\varphi)}{[0,1]} = \frac{\av{f_{\pr}(\varphi_{\pr})}{[0,1]}}{\av{(\varphi_1)_{\pr}}{[0,1]}}$. This shows that~$\Bell_{\Omega_{\pr},f_{\pr}}(y_1,y_2) = \frac{\Bell_{\Omega,f}(x_1,x_2)}{x_1}$.  

We see that if the main theorem holds true, then its analogue for the domain ``with one cusp'' (as described above) is also valid. Moreover, since the said projective transform does not change the ordering of the fixed boundary, the assertion ``the monotonic rearrangement of a function belonging to~$\Class_{\Omega}$ is also in~$\Class_{\Omega}$'' also holds true for the domains ``with one cusp''\footnote{Surely, the same method allows to cope with domains with two ``smooth'' cusps.}.

\section{Proofs of geometric statements}\label{SPG}
Lemma~\ref{HereditaryProperty} can be easily deduced from a similar one-dimensional fact.
\begin{Fact}\label{HereditaryOnLines}
Let~$F\colon [a,b] \to \mathbb{R}$ be a concave function. Suppose that~$[c,d] \subset [a,b]$. Let~$G\colon [c,d] \to \mathbb{R}$ be a concave function such that~$G(c) = F(c)$ if~$c \ne a$\textup,~$G(d) = F(d)$ if~$d \ne b$\textup, and~$G \leq F$ on~$[c,d]$. Then the function~$H$ given by the formula
\begin{equation*}
H(x) = 
\begin{cases}
G(x), \quad & x \in [c,d];\\
F(x), \quad & x \in [a,b] \setminus [c,d]
\end{cases}
\end{equation*}
is concave.
\end{Fact}
\paragraph{Proof of Lemma~\ref{HereditaryProperty}.}
We claim that the function
\begin{equation*}
G(x) = 
\begin{cases}
\BG_{w,f}(x), \quad & x\notin \tilde{w};\\
\min\big(\BG_{w,f}(x),\BG_{\tilde{w},\tilde{f}}(x)\big), \quad & x \in \tilde{w}
\end{cases}
\end{equation*}
is locally concave on~$w$. If this assertion is proved, then the lemma follows by the minimality of~$\BG_{w,f}$, because~$G \geq f$ on the fixed boundary of~$w$. To prove that~$G$ is locally concave, it suffices to prove that its restriction to each segment~$\ell \subset w$ is concave.  The set~$\ell \cap \tilde{w}$ is a closed convex subset of~$\ell$, therefore, it is either a closed interval, or a point, or an empty set. In the latter two cases, there is nothing to prove, because in such a case~$G\big|_{\ell} = \BG_{w,f}\big|_{\ell}$. The first case is Fact~\ref{HereditaryOnLines} exactly.  \qed

\begin{Rem}
The assertion of the lemma is not true if~$\tilde{w}$ is not strictly convex\textup, but only convex. 
\end{Rem}

The following fact is useful for the proof of Statement~\ref{ContinuityOnFixedFoundary}.
\begin{Fact}\label{SectionContinuity}
Suppose that a strictly convex closed set~$w \subset \mathbb{R}^d$ touches the hyperplane~$\{x\,|\, x_d = 0\}$ at the point~$0$. Then\textup, for every~$\delta > 0$ there exists~$\delta' > 0$ such that~$w \cap \{x\,|\, x_d = \pm \zeta\} \subset B_{\delta}(0)$ for all~$\zeta$\textup,~$0 < \zeta < \delta'$. 
\end{Fact}
\paragraph{Proof of Statement~\ref{ContinuityOnFixedFoundary}.}
By virtue of Lemma~\ref{HereditaryProperty}, we may change~$w$ for~$B_r(x) \cap w$ without changing~$\BG$; till the end of the proof the latter set is also denoted by~$w$. Consider a supporting hyperplane to~$w$ at~$x$, let it be the plane~$\{y \mid y_d = 0\}$. Suppose also that~$y_d \geq 0$ for~$y \in w$. Take any~$\eps > 0$. Let~$\delta$ be so small that~$|f(y) - f(x)| \leq \eps$ for all~$y \in \partial w \cap B_{\delta}(x)$. Now take~$\delta'$ as in Fact~\ref{SectionContinuity}. Then, every point~$y$ in~$w\cap B_{\delta'}(x)$ can be represented as
\begin{equation*}
y = \alpha_1y^1 + \alpha_2 y^2; \quad \alpha_1 + \alpha_2 = 1, \quad \alpha_j \geq 0, \quad y^j \in \partial w, \quad y^j_d = y_d.
\end{equation*}
In such a case,
\begin{equation*}
\BG(y) \geq \alpha_1\BG(y^1) + \alpha_2\BG(y^2) \geq f(x) - \eps.
\end{equation*}

Thus, it suffices to find~$\delta_1$ such that~$\BG(y) \leq f(x) + 2\eps$ for all~$y \in w \cap B_{\delta_1}(x)$ provided~$\eps$ is given. This~$\delta_1$ can be found in such a fashion: build a linear function~$L$ such that~$L(x) \leq f(x) + \eps$ and~$L(y) \geq \BG(y)$ for all~$y \in w$. By the minimality of~$\BG$, the inequality~$L(y) \geq \BG(y)$ may be verified only on~$\partial w$. 
The function~$L$ we are looking for is given by formula
\begin{equation*}
L(y)= f(x)+\eps+\big(\!\sup\limits_{z \in \partial \omega}|f(z)|+|f(x)|\big)\frac{y_d}{\delta'}.
\end{equation*}
The statement is proved.\qed

Finally, to prove Statement~\ref{ContinuityOnFreeBoundary}, we need a lemma.
\begin{Le}\label{MinimalityOnSimplex}
Let~$\mathcal{V}$ be the simplex with the vertices~$V_j$\textup,~$j =0,1,\ldots,d$\textup, in~$\mathbb{R}^d$. Suppose that~$K$ is an open convex set such that for~$j =1,2,\ldots,d$ one has~$[V_0,V_j] \cap \cl K = V_j$. Let~$G$ be a locally concave function on~$\mathcal{V} \setminus K$. If~$G(V_j) \geq 0$ for each~$j$\textup, then~$G \geq 0$ on~$\mathcal{V} \setminus K$. 
\end{Le}
\begin{proof}
We prove the lemma by induction in~$d$.  There is nothing to prove for~$d=1$. Let us prove the assertion of the lemma assuming that it is valid for the dimension~$d-1$. 

Assume the contrary, let~$x \in \mathcal{V} \setminus K$ be such that~$G(x) < 0$. There exists a hyperplane~$L$ that passes through~$x$ and does not cross~$K$. Its intersection with~$\mathcal{V}$ is a convex set. The boundary of this intersection belongs to the set~$\partial \mathcal{V} \setminus K$. Therefore, there exists a point~$y$ in~$\partial \mathcal{V} \setminus K$ such that~$G(y) < 0$. But the set~$\partial \mathcal{V} \setminus K$ is a union of~$(d-1)$-dimensional domains of the same type as we are working with. Therefore, the inequality~$G(y) < 0$ contradicts the induction hypothesis.
\end{proof}
\paragraph{Proof of Statement~\ref{ContinuityOnFreeBoundary}.} 
Let~$x_0 = 0$. Choose a supporting plane~$S$ to~$B_{r}(0) \setminus w$ at the point~$0$. We may assume that ~$S$ is~$\{x \mid x_d = 0\}$ and $B_{r}(0) \setminus w$ contains the point~$(0,0,\ldots,\delta)$ for some~$\delta > 0$. 
We note that the restriction of~$\BG$ to~$S$ is continuous  in a neighborhood of~$0$ by Fact~\ref{ContinuityInInnerPoints}; we take~$r$ so small that~$B_r(0) \cap S$ is contained in this neighborhood. 


To prove the continuity of~$\BG$ at~$0$, we have to prove two estimates for this function in some neighborhood of~$0$. We begin with an estimate from above. It is obtained by different methods for the points lying in~$w \cap \{x\mid x_d < 0\}$ and~$w \cap \{x\mid x_d > 0\}$. We consider the former case first.

Fix~$\eps > 0$. Choose~$\delta$ to be so small that~$\big|\BG(x_1,x_2,\ldots,x_{d-1},0) - \BG(0)\big| \leq \frac{\eps}{4}$ when~$|x_i| \leq \delta$,~$i = 1,2,\ldots,d-1$. Choose  a positive number~$\tilde{\delta} < \delta$ to be so small that~$\BG(\delta,0,0,\ldots,-x_d) \geq \BG(0) - \frac{\eps}{2}$ when~$0 < x_d< \tilde{\delta}$ (here we use the continuity of~$\BG$ at the point~$(\delta,0,\ldots,0)$), we also require the convex hull of~$B_{\tilde{\delta}}\big((\delta,0,\ldots,0)\big)$ and~$B_{\tilde{\delta}}\big((\frac{\delta}{2},0,\ldots,0)\big)$ to be in~$w$. So, if~$|x| < \delta'$ and~$x_d$ is negative, using local concavity, we get the estimate:
\begin{equation*}
\frac12\BG(x_1,x_2,\ldots,x_d) \leq \BG\Big(\frac{x_1 + \delta}{2},\frac{x_2}{2},\ldots,\frac{x_{d-1}}{2},0\Big) - \frac12\BG(\delta,0,\ldots,0,-x_d) \leq \BG(0) + \frac{\eps}{4} -\frac12\BG(0) + \frac{\eps}{4}.
\end{equation*}
Thus,~$\BG(x) \leq \BG(0) + \eps$.

Now we have to prove the estimate from above for the points in~$\{x\mid x_d > 0\}$. Fix~$\eps_0 > 0$ such that~$(0,0,\ldots,-\eps_0) \in B_r(0) \cap w$. Fix~$\eps >0$. First, we take~$\delta$ to be so small that the point~$(0,0,\ldots,-\eps_0)$  ``sees'' all the points in~$B_{\delta}(0) \cap w$ (i.e. the segment joining them is inside~$w$). Second, for~$y\in B_\delta(0)\cap \omega$ such that~$y_d>0$, we write an estimate
\begin{equation*}
\alpha_1 \BG(0,0,\ldots,-\eps_0) + \alpha_2\BG(y) \leq \BG(\alpha_2y_1,\alpha_2y_2,\ldots,0),\quad\hbox{where}\,\, -\alpha_1 \eps_0 + \alpha_2 y_d = 0,\,\,\, \alpha_1 + \alpha_2 = 1.
\end{equation*}
By continuity on~$S$,~$\BG(\alpha_2y_1,\alpha_2y_2,\ldots,0) \leq \BG(0) + \eps$ for~$\delta$ small enough, and~$y_d$ is sufficiently small. We are done.

\begin{figure}[h!]
\includegraphics[height=9.5cm]{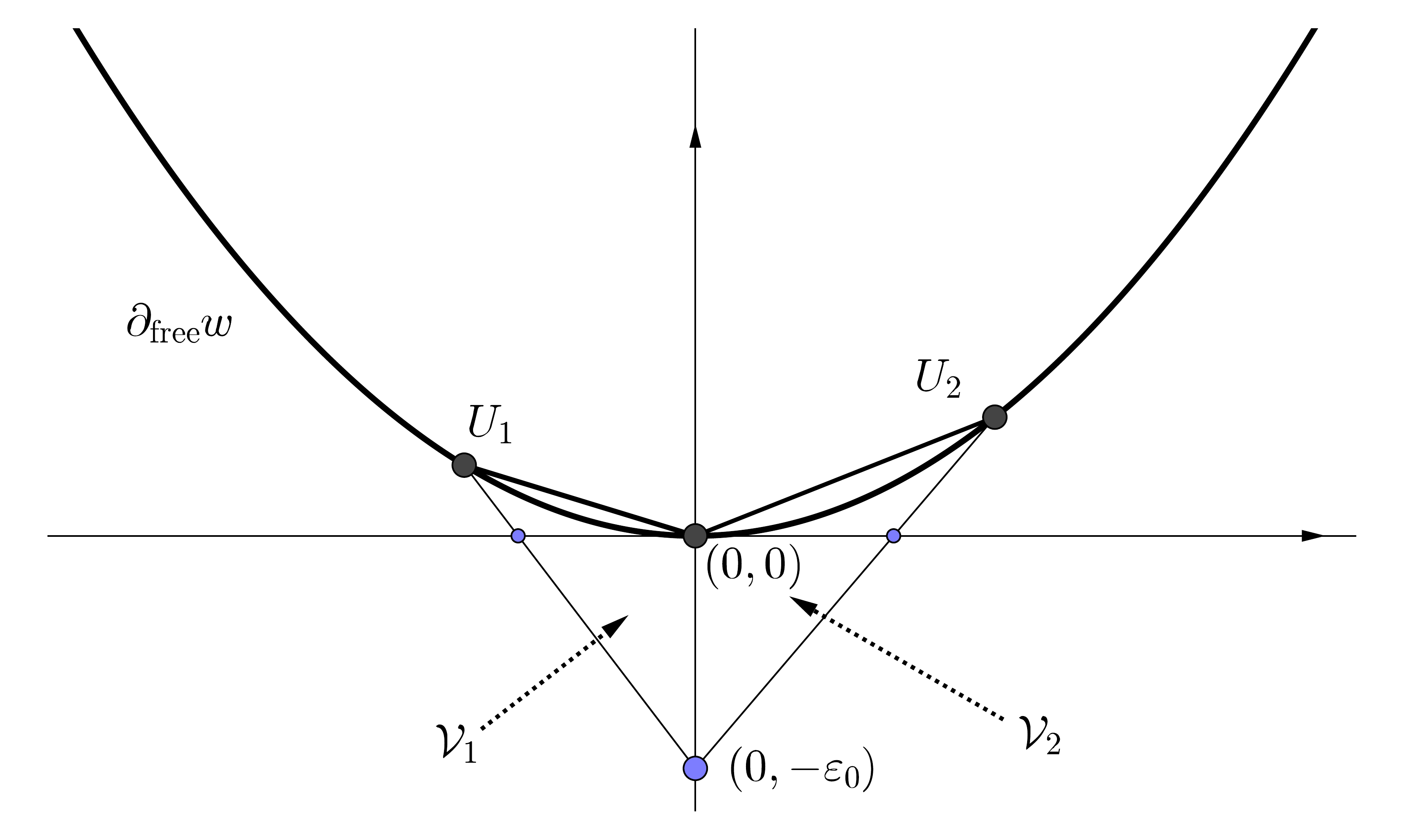}
\caption{Simplices for~$d=2$.}
\label{fig:simplices}
\end{figure}

To provide an estimate from below, we are going to cover the set~$w \cap B_{\delta}(0)$ by a finite number of simplices~$\mathcal{V}_{k}$,~$k = 1,2,\ldots,d$. Let~$\mathcal{W}$ be a simplex in~$S \cap B_{\delta}(0)$ containing~$0$. Consider the rays from~$(0,0,\ldots,-\eps_0)$ to all the vertices of~$\mathcal{W}$. Let them intersect the free boundary at the points~$U_k$,~$k=1,2,\ldots,d$. The simplex~$\mathcal{V}_k$ is a convex hull of~$(0,0,\ldots,-\eps_0)$, all the~$U_j$ except for~$U_k$, and~$0$. These simplices have two properties (see Figure~\ref{fig:simplices} also). First, their union covers~$B_{\delta'}(0) \cap w$ for every~$\delta'$ small enough. Second, each such simplex satisfies the hypothesis of Lemma~\ref{MinimalityOnSimplex}, if we take~$K$ to be~$B_{\delta}(0) \setminus w$ and~$(0,0,\ldots,-\eps_0)$ to be the~$0$-th vertex of the simplex. Therefore, for each point~$y \in B_{\delta'} \cap w$ we have an estimate
\begin{equation*}
\BG(y) \geq L_k(y), 
\end{equation*}
where~$L_k$ is a linear function interpolating~$\BG$ on the vertices of the~$\mathcal{V}_k$ the point~$y$ belongs to. This is the estimate from below we are looking for.\qed

\section{Martingales on cheese domains}\label{SCD}
We begin with two definitions.
\begin{Def}\label{MartingalesOfFiniteTime}
Let~$w$ be a subset of~$\mathbb{R}^d$. We say that an~$w$-martingale~$M$ is of finite time if for almost all~$\sigma \in \mathfrak{S}$ one has~$M_n(\sigma) \in \FixedBoundary w$ for some~$n$. 
\end{Def}
\begin{Def}\label{MartingaleConnectedDomains}
We say that a set~$w \subset \mathbb{R}^d$ is martingale connected\textup, if for each~$x \in w$ there exists a finite time martingale starting at~$x$.
\end{Def}
Our aim is to provide easy sufficient geometric conditions for a set that give the minimal principle (as Lemma~\ref{MaximalPrinciple}) and the duality theorem we discuss in Subsection~\ref{s23}. 
\begin{Def}\label{CheeseDomains}
Let~$\Omega_j$ be strictly convex compact subsets of~$\mathbb{R}^d$ with non-empty interior. Suppose that~$\Omega_j \subset \interior\Omega_0$\textup,~$j > 0$\textup, and~$\Omega_j$ are mutually uniformly separated for~$j > 0$ \textup(in particular\textup, there is only a finite number of them\textup). Then\textup,~$w = \Omega_0 \setminus \big( \cup_j \interior \Omega_j\big)$ is called a cheese domain.
\end{Def}
\begin{Th}\label{CheeseDomainsAreMartingaleConnected}
Cheese domains are martingale connected.
\end{Th}
\begin{proof}
For each point~$x$ on the free boundary of~$w$ choose an open segment~$\ell(x)$ such that~$\ell(x) \subset w$ and the endpoints of~$\ell(x)$ lie on the boundary (fixed or free) of~$w$. Such segments have a nice property: the distance from~$x$ to each of the endpoints of~$\ell(x)$ is uniformly bounded from below, because it is not less than~$\min\limits_{i,j \ne 0} \big(\min(\dist(\Omega_j,\mathbb{R}^d \setminus \Omega_0), \dist(\Omega_i,\Omega_j))\big)$.

Let~$z \in w$, we are going to construct a finite time~$w$-martingale~$M$ starting at~$z$. So, let~$M_0 =z$. We can easily split~$z$ into a convex combination~$z = \alpha_1m_1 + \alpha_2m_2$ such that~$[m_1,m_2] \subset w$ and~$m_1$ and~$m_2$ lie on the boundary of~$w$. We take~$M_1$ equal~$m_1$ with probability~$\alpha_1$ and~$m_2$ with probability~$\alpha_2$. Each of the points~$m_1$ and~$m_2$ can be splitted into two points~$m_{11}$ and~$m_{12}$ and~$m_{21}$ and~$m_{22}$ along the segments~$\ell(m_1)$ and~$\ell(m_2)$ correspondingly, if it lies on the free boundary. In the opposite case, the martingle's trajectory has reached the fixed boundary, so it does not walk any more. Then, the random variable~$M_2$ equals one of these four (or less) points with the corresponding probabilities. To pass to~$M_3$, we divide each of the points~$m_{ij}$ that lie on the free boundary, and so on. 

We have to prove that the process described above stops almost surely. We consider the values~$\E G(M_n)$, where~$G(x) = -|x|^2$. By Lemma~\ref{BellmanInduction}, the function~$n \mapsto \E G(M_n)$ is non-increasing. Moreover, the function~$G$ is uniformly strictly concave, so, we can say more. If~$x_+$ and~$x_-$ are in~$\mathbb{R}^d$ and~$\alpha_+$ and~$\alpha_-$ are separated from~$0$ and~$1$, then
\begin{equation*}
G(\alpha_+x_+ + \alpha_-x_-) \geq \alpha_+G(x_+) + \alpha_-G(x_-) + \tilde{c}|x_+ - x_-|^2
\end{equation*}    
for some small~$\tilde{c}$. Therefore, if~$M_n$ is in~$\FreeBoundary w$ with probability~$p_n$, then
\begin{equation*}
\E G(M_{n+1}) \leq \E G(M_n) - cp_n,
\end{equation*}
because the lengths of the segments~$\ell(x)$ are bounded from below. The value~$\E G(M_n)$ is bounded from below (because~$G$ is bounded from below on~$\Omega_0$), so~$p_n \to 0$, which means that~$M_n$ stops eventually a.s.
\end{proof}
\begin{Le}[{\bf Minimal principle for cheese domains}]\label{MaximalPrincipleForCheeseDomains}
Let~$w$ be a cheese domain. Suppose that the function~$f\colon \FixedBoundary w \to \mathbb{R}$ is continuous. Then\textup, 
\begin{equation*}
\min_{x\in w}\BG_{w,f}(x) = \min_{x\in \FixedBoundary w} f(x).
\end{equation*}
\end{Le}
\begin{proof}
We note that cheese domains satisfy the assumptions of Propositions~\ref{ContinuityOnFixedFoundary} and~\ref{ContinuityOnFreeBoundary} that assert that~$\BG$ is continuous on~$w$. As in the proof of Lemma~\ref{MaximalPrinciple}, we show that~$\BG(z) \geq \min_{x \in \FixedBoundary w} \BG(x)$ for all~$z \in w$. Let~$M^z$ be a finite time~$w$-martingale starting at~$z$, the existence of such a martingale is provided by Theorem~\ref{CheeseDomainsAreMartingaleConnected}. By Lemma~\ref{BellmanInduction},~$\BG(z) \geq \E \BG(M_n)$. The latter value tends to~$\E f(M_{\infty}^z)$ by the Lebesgue dominated convergence theorem (we recall that~$\BG$ is a continuous function on a compact set), which is not smaller than~$\min_{x \in \FixedBoundary w} f(x)$.
\end{proof}
The reader can easily construct sets~$w$ for which the minimal principle does not hold. For example, if~$\Omega$ is as described in the beginning of Section~\ref{s1}, but does not satisfy condition~\eqref{ThirdCondition} (and the set~$\Omega_1$ is not bounded), then the minimal principle does not hold for it. Moreover, there exist compact sets~$w$ for which the minimal principle does not hold. However, all the compact sets the authors managed to construct have fractal boundary. 
\begin{Conj}\label{Manifold}
Suppose that~$w$ is a bounded closed connected set with non-empty interior\textup, whose boundary is a real-analytic submanifold of~$\mathbb{R}^d$. Then\textup,~$w$ is martingale connected\textup.
\end{Conj}

We end this section with an analog of Theorem~\ref{GeometricMartingale}.

\begin{Th}\label{GeometricMartingaleCheeseDomain}
Let~$w$ be a cheese domain. Suppose that~$f\colon \FixedBoundary w \to \mathbb{R}$ is continuous. Then\textup,~$\BG_{w,f} = \BM_{w,f}$.
\end{Th}
\begin{proof}
The proof is the same as for Theorem~\ref{GeometricMartingale}, one only has to recall Proposition~\ref{ContinuityOnFixedFoundary} for the continuity of the function~$\BG$ at the fixed boundary and replace Lemma~\ref{MaximalPrinciple} by Lemma~\ref{MaximalPrincipleForCheeseDomains} in the proof of the analog of Lemma~\ref{BellmanInductionLimit}.
\end{proof}




Dmitriy M. Stolyarov
\medskip

P. L. Chebyshev Research Laboratory, St. Petersburg State University,

St. Petersburg Department of Steklov Mathematical Institute, Russian Academy of Sciences;
\medskip

dms at pdmi dot ras dot ru.
\bigskip

Pavel B. Zatitskiy 
\medskip

P. L. Chebyshev Research Laboratory, St. Petersburg State University,

St. Petersburg Department of Steklov Mathematical Institute, Russian Academy of Sciences;
\medskip

paxa239 at yandex dot ru.

\end{document}